\documentclass[11pt]{amsart}
\allowdisplaybreaks[1]

\usepackage{a4wide}
\usepackage{amssymb}
\usepackage{amsmath}
\usepackage{graphicx}

\usepackage{graphicx}

\usepackage{epstopdf}

\input xypic

\newtheorem{theorem}{Theorem}[section]
\newtheorem{lemma}[theorem]{Lemma}
\newtheorem{prop}[theorem]{Proposition}
\newtheorem{coro}[theorem]{Corollary}

\theoremstyle{definition}

\newtheorem{remark}{Remark}

\newcommand{\RR}{\mathbb{R}}

\newcommand{\CC}{\mathbb{C}}

\newcommand{\TT}{\mathbb{T}}

\newcommand{\ZZ}{\mathbb{Z}}

\newcommand{\om}{\omega}

\def\a{\alpha}
\def\o{\omega}
\def\r{\rangle}
\def\l{\langle}

\newcommand{\Hm}[1]{\leavevmode{\marginpar{\tiny%
$\hbox to 0mm{\hspace*{-0.5mm}$\leftarrow$\hss}%
\vcenter{\vrule depth 0.1mm height 0.1mm width \the\marginparwidth}%
\hbox to
0mm{\hss$\rightarrow$\hspace*{-0.5mm}}$\\\relax\raggedright #1}}}

\title[Cubature formulae for simple Lie groups]
{Cubature formulae for orthogonal polynomials\\ in terms of\\ 
elements of finite order of compact simple Lie groups}

\author{Robert V.~Moody}
\address{Department of Mathematics and Statistics,
University of Victoria, \newline
\hspace*{12pt}Victoria, BC, V8W3P4, Canada}
\email{rmoody@uvic.ca}

\author{Ji\v{r}\'\i\ Patera}
\address{CRM,
Universit\'e de Montr\'eal, Montr\'eal}
\email{patera@crm.umontreal.ca}

\date{\today}
\thanks{}

\begin{document}
\begin{abstract}
The paper contains a generalization of known properties of Chebyshev polynomials of the second kind in one variable  to polynomials of $n$ variables based on the root lattices of compact simple Lie groups $G$ of any type and of rank $n$. The results, inspired by work of H. Li and Y. Xu where they derived cubature formulae from $A$-type lattices, yield Gaussian cubature formulae for each simple Lie group $G$ based on interpolation points that arise from regular elements of finite order in $G$. The polynomials arise from the irreducible characters of $G$ and the interpolation points as common zeros of certain finite subsets of these characters. The consistent use of Lie theoretical methods reveals the central ideas clearly and allows for a simple uniform development of the subject. Furthermore it points to  genuine and perhaps far reaching Lie theoretical connections. 

\end{abstract}

\maketitle

\section{Introduction}\label{intro}

During most of the century and half long history of Chebyshev polynomials, only polynomials of one variable were studied \cite{Rivlin}. In recent years a considerable overlap of the subject can be found in the emerging and practically important field of cosine and sine transforms \cite{Rao,Strang}. In the absence of additional constraints, truly higher dimensional generalizations of Chebyshev polynomials are hidden in the vast number of possibilities of defining orthogonal polynomials of more than one variable \cite{suetin, Mac2}. In $2D$ \cite{Koornwinder1-4} such constraints were provided by requiring that the polynomials of two variables be simultaneous eigenvectors of two differential operators. 

Our motivation here comes from two directions: 
\begin{itemize}
\item[(i)] The method of constructing Chebyshev-like polynomials for any simple Lie group \cite{NPT}, and particularly from the recognition of the basic role played by
the group characters. 
\item[(ii)] The work of H.~Li and Y.~Xu in which they derived cubature formulae based on the symmetries of $A$-type lattices \cite{LX}. 
\end{itemize}
From \cite{NPT} we understood that a general formulation, uniform over all the types and ranks of simple Lie algebras should be possible, and from \cite{LX} we saw what possibilities, beyond the construction of the polynomials, might be achievable in a general formulation guided by the theory of compact simple Lie groups.

In this paper the characters of the irreducible representations of a compact simply-connected simple Lie group $G$ of rank $n$ play the role of the Chebyshev polynomials (of the second kind) and elements of finite order \cite{MP84} in $G$ give rise to interpolation points through which we arrive at Gaussian cubature formulas. The ring generated by the characters of  $G$ has a $\ZZ$-basis consisting of the irreducible characters and it is a polynomial ring in terms of the $n$ irreducible characters of the fundamental representations. These fundamental characters then serve as new variables for functions defined on a bounded domain $\Omega \subset \RR^n$. This domain $\Omega$ is derived from the fundamental domain of the affine Weyl group of $G$ and the kernel function used is the absolute value of the denominator
of Weyl's character formula \cite{BS}.  

There are two technical ingredients in the paper, which are indispensable for the uniformity of our approach to simple Lie groups of all types and thus for generality of our conclusions.
\smallskip
\newline
(i) The `natural' grading of polynomials by their total degree is replaced by new $m$-degree grading. It is based on a set of Lie theoretical invariants, called the marks, which are unique for simple Lie group of each type. The two gradings coincide only in the case of $A_n$. 
\smallskip
\newline
(ii) The set of interpolation points is uniformly specified for each simple Lie group as a finite set of lattice points, characterising all conjugacy classes of elements of certain order in the underlying Lie group.
\smallskip

 The cubature formula (see Cor.~\ref{cubature}) is then a formula that equates an integral of a general polynomial $P$, say of $m$-degree $\leq 2M+1$, to weighted finite sums of the values of $P$ sampled at the interpolation points which are common zeros of the polynomials of $m$-degree $M+1$. The cubature is {\em Gaussian} \cite{Xu} in the sense that the number of the interpolation points coincides with the dimension of the space of polynomials of $m$-degree $\leq M$. It turns out that the interpolation points are the elements whose adjoint order divides $M+h$ where $h$ is the Coxeter number.

In more detail, having fixed a specific $m$-degree $M$, we are interested in the properties of the set of all polynomials of $m$-degree at most $M$. There are three main results:

\begin{itemize}
\item[{1)}] Thm.~\ref{0points}: The interpolation points are common zeros for the set of polynomials of $m$-degree $M+1$. These interpolation points are directly related to regular elements of finite adjoint order $M+h$ in $G$, where $h$ is the Coxeter number of $G$, and their number is precisely the dimension of the space of polynomials of $m$-degree less than or equal to $M$. 

\item[{2)}] Cor.~\ref{cubature}: There is a cubature formula that equates weighted integrals 
$\int_\Omega f \, K^{1/2} $ of each polynomial $f$ of $m$-degree $\leq2M+1$ with 
$K$-weighted linear combinations of its values at the interpolation points. 

\item[{3)}] Prop.~\ref{approxThm}: The expansions of functions via irreducible characters (or rather these
characters interpreted as polynomials) using interpolation points yield the best approximations in the 
$K^{1/2}$-weighted $L^2$-norm on the functions on $\Omega$ .
 \end{itemize}
 
\smallskip
These results and their proofs are natural and completely uniform within the framework of the character theory of simple Lie groups and their elements of finite order. In fact it is this naturalness and perfection of fit that suggests that there are a deeper Lie theoretical implications to all of this that are still to be discovered. 
The potential role of the finite reflection groups in the theory of orthogonal polynomials has been recognized already in \cite{DX}, although only a limited use of the groups is made there, the objects of interest being group invariant differential and difference-differential operators.  The present paper can be understood as a new contribution to the fulfilment of that potential, much closer to the properties of the simple Lie groups which give rise to the finite reflection groups. 

Our discussion requires a certain familiarity with root and weight lattices, their corresponding co-root and co-weight lattices, the Weyl group and its affine extension, and particularly the structure of the natural fundamental domain  for the affine Weyl group. 
We use the \S\ref{Lie groups}, which sets up the notation, to briefly review the key points of this theory, while \S\ref{Winvariance} and \S\ref{EFOs} contain preparatory extensions of the standard theory. Consideration of the subject of the paper really starts in \S\ref{start}.

\section{Simple compact Lie groups} \label{Lie groups}

This section is a brief review of the material that we need for this paper and establishes
the notation that we shall be using. The main facts about simple Lie groups and their
representations are classical and can be found in many places. One source
that uses the same notation as in this paper is \cite{KMPS}, Vol.1. Material on the fundamental domain, and in particular information about the stabilizers
of its various points, can be found in Ch.V of \cite{B} and material on both the fundamental domain and the elements of finite order can be found in \cite{MP87}.

Let $G$ be a simply connected simple Lie group with Lie algebra $\mathfrak{g}$. We let
$c_G$ denote the order of its centre and let $\TT$ be a maximal torus of $G$. We let 
$i\mathfrak{t}$ be the Lie algebra of $\TT$, so that we have the exact sequence
$$
 0 \longrightarrow \check Q \longrightarrow \mathfrak t 
\stackrel{\exp(2 \pi i (\cdot))}{\longrightarrow} \TT \longrightarrow 1
$$
via the exponential map. Using $\mathfrak{t}$ instead of the actual Lie algebra $i\mathfrak{t}$ of $\TT$ has the advantages that the Killing form $(\cdot\mid\cdot)$ of $G$ restricted to $\mathfrak{t}$ is positive definite and $e^{2\pi i(\cdot\mid\cdot)}$ is perfect for the Fourier analysis to follow. Let $n:= \dim_\RR\mathfrak t$, the {\em rank} of $G$.

The kernel of $e^{2\pi i(\cdot)}$ on $\mathfrak{t}$ is the {\em co-root lattice} $\check Q$, and $\mathfrak{t}/\check Q\simeq\TT$ naturally expresses $\TT$ as a real space factored by a lattice. We denote by $\mathfrak t^*$ the dual space of $\mathfrak t$ and let
$\langle \cdot\,,\cdot \rangle$ be the natural pairing of $\mathfrak t^*$ and $\mathfrak t$.

Let $\theta_\TT$ be the Haar measure on $\TT$ that gives it volume equal to $1$. In practice
we most often write integration over $\TT$ in the form of integration over some fundamental region $FR$
for $\check Q$ in $\mathfrak t$:
\begin{equation}\label{integralsOverT}
\int_\TT f \, d\theta_\TT = \int_{FR} f(e^{2\pi i x}) d\theta_{\mathfrak t} 
\end{equation}
where $\theta_{\mathfrak t}$ is ordinary Lebesgue measure in $\mathfrak t$ normalized
so that $FR$ has volume equal to $1$. In the sequel no fundamental domain $FR$ ever makes an appearance, but rather we work with a smaller fundamental domain $F^\circ$ of the affine Weyl group,
see below.

Given any finite dimensional complex representation $V$ of the Lie group $G$, the action
of the elements of $\TT$ on $V$ can be simultaneously diagonalized, the resulting eigenspaces in $V$ being the {\em weight spaces}. The corresponding action of $\mathfrak t$
on these weight spaces then affords elements $\lambda \in \mathfrak t^*$ for which
$x\in \mathfrak t$ acts as $\exp(2 \pi i \langle \lambda, x\rangle)$. Naturally 
$\langle \lambda, \check Q \rangle \subset \ZZ$ and the set of weights taken over all
finite dimensional representations is the subgroup $P \subset
\mathfrak t^*$ which is the $\ZZ$-dual of the co-root lattice relative to our pairing. $P$
is the {\em weight lattice} of $G$ relative to $\TT$. The adjoint representation, the 
representation of $G$ on its own Lie algebra, produces its own set of weights. Apart from the weight $0$, which is of multiplicity $n$, the remaining weights of the adjoint representation
are the {\em roots} of $G$ and they all occur with multiplicity $1$. They generate the
root lattice $Q$ of $G$, and we have $Q\subset P$ with index equal to $c_G$ defined
above. The $\ZZ$-dual of $Q$ in $\mathfrak t$ is the co-weight lattice $\check P
\supset \check Q$ with the index of $\check Q$ in $\check P$ being $c_G$, again.

Let $W$ be the Weyl group associated with $\TT$, namely the quotient of the normalizer of $\TT$ in $G$. Then $W$ acts as a finite group of isometries of $\mathfrak t$ 
relative to $(\cdot\mid\cdot)$, and $W$ stabilizes $\check Q$.  It also stabilizes $\check P$, and if we pull its action over to $\mathfrak t^*$ by
duality, then $W$ also stabilizes the weight lattice $P$ and the root lattice $Q$.

The relationships between the lattices and between the various root and weight bases and their co-equivalents described below are summarized in:
\begin{equation}
\begin{matrix}
\{\alpha_1, \dots, \alpha_n\}&\subset &Q&&\check Q&&\\
 &&\cap&\times&\cap\\
& &P&&\check P&\supset&\{\check\omega_1, \dots, \check\omega_n\}\\
&&\cap&&\cap&\\
&& \mathfrak t^*&& \mathfrak t
\end{matrix}
\end{equation}
The times symbol is meant to indicate that $Q$ and $\check P$, as well as
$\check Q$ and $P$, are in $\ZZ$-duality with each other. The indicated bases are
also in $\ZZ$-duality. 

\smallskip
We have the semi-direct product $W_{\mbox{\small{aff}}}= W \ltimes\check Q$ acting on $\mathfrak{t}$, with $\check Q$ acting as translations and the Weyl group as point symmetries. $W_{\mbox{\small{aff}}}$ is called the {\em affine Weyl group}.

\begin{figure}
 \hspace{-5cm}
\parbox{.6\linewidth}{\setlength{\unitlength}{2pt}
\def\kr{\circle{4}}
\def\cr{\circle*{4}}
\def\Kr{\circle{4}}
\thicklines
\begin{picture}(210,110)

\put( 0,95){\makebox(0,0){${A_n}$}}
\put(10,95){\kr}
\put(20,95){\kr}
\put(30,95){\kr}\put(30,102){\Kr}\put(29.4,100.5){$\cdot$}
\put(37,95){$\dots$}
\put(50,95){\kr}
\put(12,95){\line(1,0){6}}
\put(22,95){\line(1,0){6}}
\put(32,95){\line(1,0){4}}
\put(44,95){\line(1,0){4}}
\put(11.5,96){\line(3,1){16.5}}
\put(48.5,96){\line(-3,1){16.5}}
\put(68,93){$n\geq 1$}

\put( 0,70){\makebox(0,0){${B_n}$}}
\put(10,70){\kr}
\put(20,70){\kr}\put(19,61){$\tfrac22$}\put(20,72){\line(0,1){4}}
   \put(20,78){\Kr}\put(19.4,76.5){$\cdot$}
\put(27,70){$\dots$}
\put(40,70){\kr}\put(39,61){$\tfrac22$}
\put(50,70){\cr}\put(49,61){$\tfrac22$}
\put(12,70){\line(1,0){6}}
\put(22,70){\line(1,0){4}}
\put(34,70){\line(1,0){4}}
\put(42,71){\line(1,0){6}}
\put(42,69){\line(1,0){6}}
\put(68,69){$n\geq 2$}

\put( 0,45){\makebox(0,0){${C_n}$}}
\put(10,45){\Kr}\put(9.3,43.5){$\cdot$}
\put(20,45){\cr}\put(19,36){$\tfrac22$}
\put(30,45){\cr} \put(29,36){$\tfrac22$}
\put(37,45){$\dots$}
\put(50,45){\cr}\put(49,36){$\tfrac22$}
\put(60,45){\kr}
\put(22,45){\line(1,0){6}}
\put(32,45){\line(1,0){4}}
\put(44,45){\line(1,0){4}}
\put(12,46){\line(1,0){6}}
\put(12,44){\line(1,0){6}}
\put(52,46){\line(1,0){6}}
\put(52,44){\line(1,0){6}}
\put(68,44){$n\geq 2$}

\put( 0,20){\makebox(0,0){${D_n}$}}
\put(10,20){\kr}
\put(20,20){\kr}\put(19,11){$\tfrac22$}\put(20,28){\Kr}\put(19.4,26.3){$\cdot$}
\put(40,20){\kr}\put(39,11){$\tfrac22$}
\put(50,20){\kr}\put(49,11){$\tfrac22$}
\put(60,20){\kr}
\put(50,28){\kr}
\put(12,20){\line(1,0){6}}
\put(22,20){\line(1,0){4}}
\put(27,20){$\dots$}
\put(34,20){\line(1,0){4}}
\put(42,20){\line(1,0){6}}
\put(52,20){\line(1,0){6}}
\put(20,22){\line(0,1){4}}
\put(50,22){\line(0,1){4}}
\put(68,19){$n\geq 4$}

\put(110,95){\makebox(0,0){${E_6}$}}
\put(120,95){\kr}
\put(130,95){\kr}\put(129,86.5){$\tfrac22$}\put(130,103){\Kr}\put(129.2,101.5){$\cdot$}
\put(140,95){\kr}\put(139,86.5){$\tfrac33$}\put(140,103){\kr}\put(143,102){$\tfrac22$}
\put(150,95){\kr}\put(149,86.5){$\tfrac22$}
\put(160,95){\kr}
\put(122,95){\line(1,0){6}}
\put(132,95){\line(1,0){6}}
\put(142,95){\line(1,0){6}}
\put(152,95){\line(1,0){6}}
\put(140,97){\line(0,1){4}}
\put(132,103){\line(1,0){6}}

\put(110,70){\makebox(0,0){${E_7}$}}
\put(120,70){\Kr}
\put(130,70){\kr}\put(129,61){$\tfrac22$}
\put(140,70){\kr}\put(139,61){$\tfrac33$}
\put(150,70){\kr}\put(149,61){$\tfrac44$}\put(150,78){\kr}\put(153,77){$\tfrac22$}
\put(160,70){\kr}\put(159,61){$\tfrac33$}
\put(170,70){\kr}\put(169,61){$\tfrac22$}
\put(180,70){\kr} \put(179.3,68.5){$\cdot$}
\put(122,70){\line(1,0){6}}
\put(132,70){\line(1,0){6}}
\put(142,70){\line(1,0){6}}
\put(152,70){\line(1,0){6}}
\put(162,70){\line(1,0){6}}
\put(172,70){\line(1,0){6}}
\put(150,72){\line(0,1){4}}

\put(110,45){\makebox(0,0){${E_8}$}}
\put(120,45){\Kr}\put(119.3,43.5){$\cdot$}
\put(130,45){\kr}\put(129,36){$\tfrac22$}
\put(140,45){\kr}\put(139,36){$\tfrac33$}
\put(150,45){\kr}\put(149,36){$\tfrac44$}
\put(160,45){\kr}\put(159,36){$\tfrac55$}
\put(170,45){\kr}\put(169,36){$\tfrac66$}
   \put(170,53){\kr}\put(173,52){$\tfrac33$}
\put(180,45){\kr} \put(179,36){$\tfrac44$}
\put(190,45){\kr}\put(189,36){$\tfrac22$}
\put(122,45){\line(1,0){6}}
\put(132,45){\line(1,0){6}}
\put(142,45){\line(1,0){6}}
\put(152,45){\line(1,0){6}}
\put(162,45){\line(1,0){6}}
\put(172,45){\line(1,0){6}}
\put(182,45){\line(1,0){6}}
\put(170,47){\line(0,1){4}}

\put(110,20){\makebox(0,0){${F_4}$}}
\put(120,20){\Kr}\put(119.3,18.3){$\cdot$}
\put(130,20){\kr}\put(129,11){$\tfrac22$}
\put(140,20){\kr}\put(139,11){$\tfrac34$}
\put(150,20){\cr}\put(149,11){$\tfrac43$}
\put(160,20){\cr}\put(159,11){$\tfrac22$}
\put(122,20){\line(1,0){6}}
\put(132,20){\line(1,0){6}}
\put(142,21){\line(1,0){6}}
\put(142,19){\line(1,0){6}}
\put(152,20){\line(1,0){6}}

\put(180,20){\makebox(0,0){${G_2}$}}
\put(190,20){\kr}\put(189,11){$\tfrac23$}\put(190,28){\Kr}\put(189.3,26.5){$\cdot$}
\put(200,20){\cr}\put(199,11){$\tfrac32$}
\put(190,22){\line(0,1){4}}
\put(192,20){\line(1,0){6}}
\put(190,21.5){\line(1,0){10}}
\put(190,18){\line(1,0){10}}

\end{picture}
}
\caption{The usual Coxeter-Dynkin diagrams, which are simply a visually efficient way to encode the Cartan matrices, are shown, along with information giving the corresponding marks and co-marks. The circular nodes, ignoring those with a dot in them, stand for the simple roots with the convention that open (resp. filled) circles indicate long (resp. short) roots. 
When there are both types of nodes, interchanging open and filled results in the co-diagram. Thus types $B$ and $C$
interchange and types $F$ and $G$ end up as the same type but with the numbering
of the roots permuted. The dotted node stands
for the root $\alpha_0$ (which is linearly dependent on the simple roots with the negatives
of the marks as its coefficients). Under the duality operation of roots to co-roots, it
passes from being the lowest long root to being the lowest short co-root, a fact that
we do not explicitly have to use here. \newline
The links between roots occur only when the roots are
not orthogonal to one another. The marks and co-marks are shown as the fractions $m/\check m$ attached to the corresponding nodes of the diagram. When both are equal to one, they are not shown. We refer the reader to \cite{BMP} for more details.
The numbering of simple roots goes from the left to right, the node above the main line carrying the highest value. Dotted cicle has number 0.}\label{diagrams}

\end{figure}
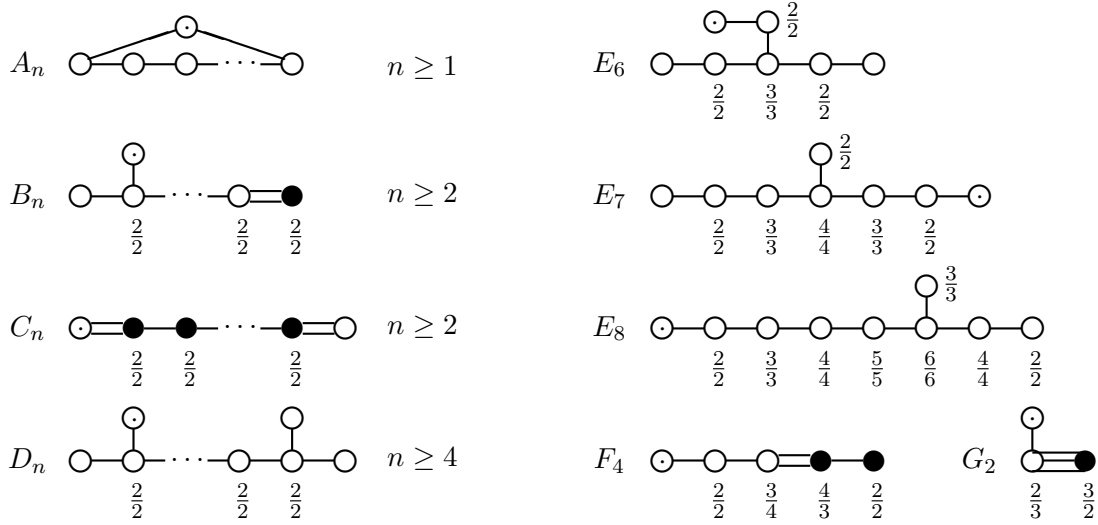

A fundamental region $F$ for $\mathfrak{t}$ under the action of $W_{\mbox{\small{aff}}}$ can be given as follows.  
Let $\Delta$ be the set of all the roots and let $\Pi=\{\alpha_1,\dots,\alpha_n\}$ be a simple system of roots for $\Delta$. Let $\Delta_+$ (resp. $\Delta_-$) denote the corresponding set of positive (resp. negative) roots, and let
\begin{equation}\label{marksDefined}
\alpha_0=-(m_1\alpha_1+\cdots+m_n\alpha_n)
\end{equation}
be the lowest root\footnote{In what follows, the use of the lowest root rather than the highest root, which might seem more natural, is in keeping with notation from the theory of affine root systems.} of $\Delta$. The positive integers $m_1,\dots,m_n$ are called the {\it marks} of $G$. The marks and co-marks are shown on Fig.~\ref{diagrams}. They are independent of the choice of $\TT$ and the choice of $\Pi$ (but depend on the choice of ordering of  $\alpha_1,\dots,\alpha_n$). It is convenient to define
one more mark, $m_0:=1$, so one has $\sum_{j=0}^n m_j\alpha_j =0$. One then
has the {\em Coxeter number} of $G$
 $$
 h:= \sum_{j=0}^n m_j \, .
 $$ 
 This constant appears in many places and in many guises in the theory, and plays an important role in what follows.

Now we define 
\begin{equation}
F:=\{x\in\mathfrak{t}\mid \l\alpha_j,x\r\geq0\text{ for all }j=1,\dots,n,\ 
\l -\alpha_0,x\r\leq1\}\,.
\end{equation}
This is a simplex in $\mathfrak{t}$.
The Weyl group $W$ contains the so-called {\em simple} reflections $r_1,\dots,r_n$ in the walls 
\begin{equation}\label{hyperplane}
H_j:=\{x\in\mathfrak{t}\mid\l\alpha_j,x\r=0;\ j=1,\dots,n\}
\end{equation}
and $W$ is generated by these reflections. The {\em length} function on $W$
is the homomorphism
$$
 l : W \longrightarrow \{\pm1\}, \quad r_j \mapsto -1 \quad\mbox{for all} \; j = 1, \dots, n \,.
 $$
There is a unique element
$w_{\mbox{\small{opp}}} \in W$ which maps $\Delta_+$ into $\Delta_-$. It is an involution
and is minimally represented by the product of exactly $|\Delta_+|$ of the simple
reflections.

\begin{figure}
\centering
\includegraphics[width=12cm]{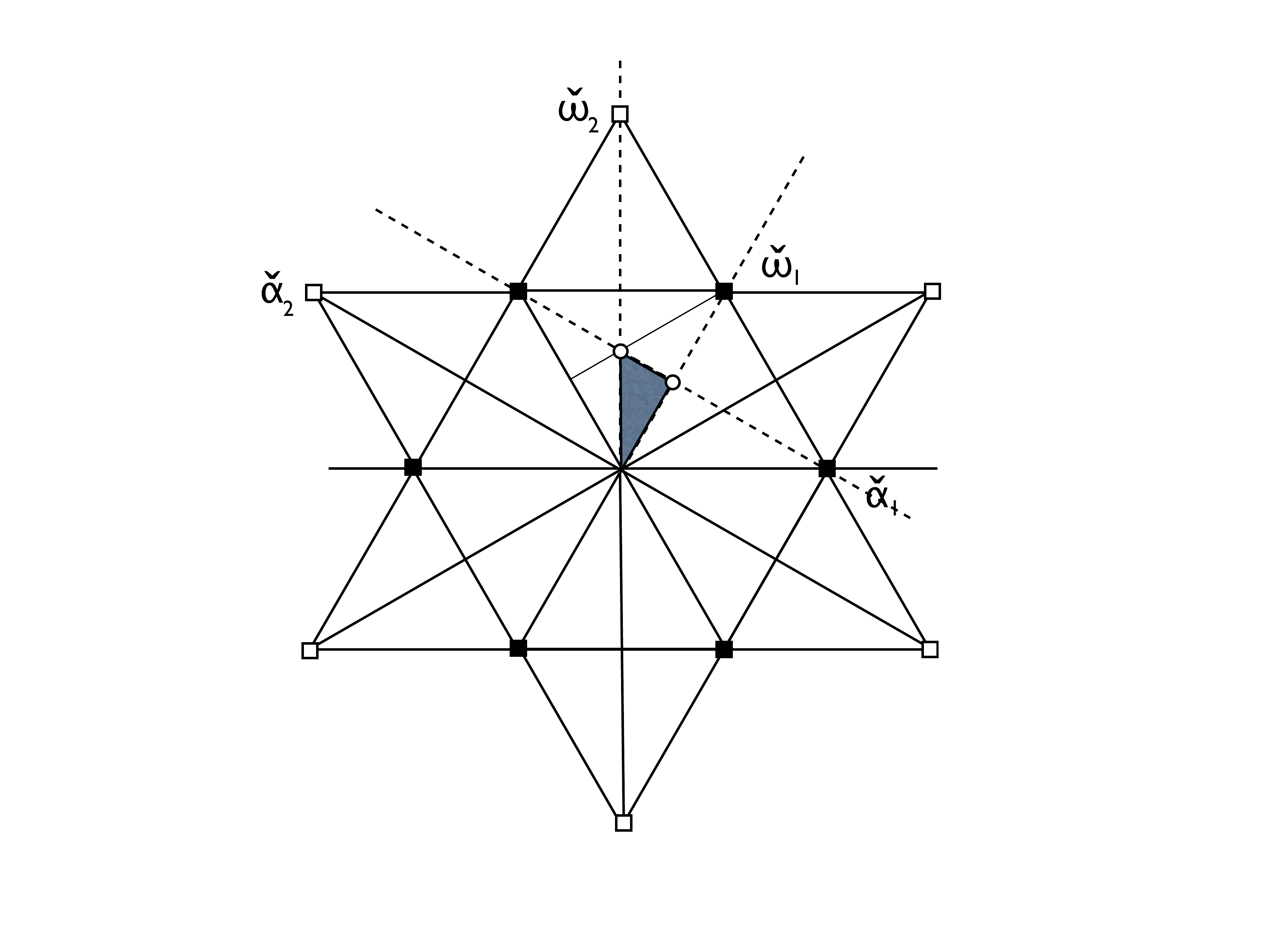}
\caption{A schematic view of the co-root system of $G_2$. The shaded triangle is the
fundamental region $F$. The dotted lines are the mirrors which define its boundaries, the
reflections in which generate the affine Weyl group. The action of the affine Weyl group 
on $F$ tiles the plane.  A few tiles of this tiling are shown. Filled (resp, open) squares are the short (resp. long) co-roots of $G_2$. }\label{G2figure}
\label{g2}
\end{figure}

$W_{\mbox{\small{aff}}}$ is also a reflection group and is obtained by adding to $r_1, \dots, r_n$ the additional generator $r_0$ which is reflection in the wall 
\begin{equation}\label{0-hyperplane}
H_0:=\{x\in\mathfrak{t}\mid\l -\alpha_0,x\r=1\}\,.
\end{equation}
There is a similar length function for $W_{\mbox{\small{aff}}}$
The action of  $W_{\mbox{\small{aff}}}$ on $F$ tiles the entire space $\mathfrak t$ with copies of 
$F$, and in this way $F$ serves as a fundamental domain for it. See Fig.~\ref{G2figure} for
an example.

The way in which $F$ is the fundamental region is rather beautiful:
\begin{itemize}
\item
Each  $W_{\mbox{\small{aff}}}$ orbit in $\mathfrak{t}$ has a  unique element in $F$.
\item
For $x\in F$ the stabilizer of $x$ in $W_{\mbox{\small{aff}}}$ is the subgroup of
$W_{\mbox{\small{aff}}}$ generated by the reflections $r_j$, \ $j=0,1,\dots,n$, for which $x$ is on the wall $H_j$. In particular, all $x\in F^\circ$ have trivial stabilizer in $W_{\mbox{\small{aff}}}$.
\end{itemize}

The {\em fundamental co-weights} are the elements $\check\o_k \in \mathfrak t$
dual to the simple roots: $\l \alpha_j,\check\o_k\r = \delta_{jk}$. They form a
$\ZZ$-basis of $\check P$. In terms of them we can write $F$ as the convex hull of
$$
\left\{0,\,\frac{\check\omega_1}{m_1},\,\frac{\check\omega_2}{m_2},\,\cdots,\,
\frac{\check\omega_n}{m_n}\right\}\,.
$$

Although each element of $W$ and each element of $W_{\mbox{\small{aff}}}$ is a product of reflections, it is not uniquely so. However, the parity of the number of
reflections is well-defined and we use the symbol $(-1)^{l(w)}$ to indicate this parity,  it being $\pm 1$ according to evenness or oddness.  

The Cartan matrix of $G$ is the integer $n\times n$ matrix 
$$
A=(A_{ij})
  =\left(\frac{2(\alpha_i\mid\alpha_j)}{(\alpha_j\mid\alpha_j)}\right)\,,\qquad
  1\leq i,j\leq n\,.
$$
It is unique up to the numbering of simple roots. The Cartan matrices classify the compact simply-connected  simple Lie groups into the well-known $A$, $B$, $C$, $D$ series and the five exceptional groups $E_6$, $E_7$, $E_8$, $F_4$, and $G_2$.

The simple {\em co-roots} $\{\check\alpha_1,\dots,\check\alpha_n\}\subset \check Q$ corresponding to $\{\alpha_1,\dots,\alpha_n\}\subset Q$ are
 defined by 
$$
\l\alpha_i,\check\alpha_j\r=A_{ij}
    =\frac{2(\alpha_i\mid\alpha_j)}{(\alpha_j\mid\alpha_j)}\qquad
    \text{for all}\quad i,j\,.
$$
They form a $\ZZ$-basis of $\check Q$ and their $W$-translates in $\check Q$ form the set of co-roots $\check\Delta$. Actually $\check\Delta$ and the simple co-roots are a root system and simple roots for a simply-connected group  $\check G$ whose Cartan matrix is $A^T$. We do not need this group directly in what follows but we occasionally use information about co-objects that we know is true from the fact that they have such an interpretation.

Dual to the co-roots we have the {\em fundamental weights} $\omega_j \in P$
defined by $\l \omega_j, \check \alpha_k\r = \delta_{jk}$. The fundamental
weights form a $\ZZ$-basis of $P$.

Each finite dimensional irreducible representation $L$ of $G$ has a unique one dimensional weight space $L^\lambda$, \ $\lambda\in P$, with the property that all other weights of $L$ are of the form $\lambda-\beta$, where $\beta$ is a sum of positive roots. We have
$$
\lambda\in P^+
:=\{\mu\in\mathfrak t^*\mid\l\mu,\check\alpha_i\r\in\ZZ^{\geq0},\ i=1,\dots,n\} \,.
$$
Here $P^+\subset P$ is the set of {\it dominant weights}. The dominant weight $\lambda$ 
is called the {\it highest weight} of $L$. If we designate $L$ now by $L(\lambda)$, then the correspondence
$$
\lambda\in P^+\longleftrightarrow L(\lambda)
$$
classifies all the irreducible finite dimensional representations of $G$ up to isomorphism.

We shall also need the set of {\em strictly dominant weights}
$$
P^{++}:=\{\mu\in\mathfrak t^*\mid\l\mu,\check\alpha_i\r\in\ZZ^{>0},\ i=1,\dots,n\}\,.
$$
The `simplest' strictly dominant weight is $\rho$ defined by $\l\rho,\check\alpha_i\r=1$ for $i=1,\dots,n$. The element  $\rho$ plays an important role in what follows. We also know that
\begin{equation}\label{rho}
\rho=\tfrac12\sum_{\alpha\in\Delta_+}\alpha=\sum_{k=1}^n\omega_k\,.
\end{equation}

The character of a finite dimensional representation $L$ of $G$ is the mapping
$$
g\ \mapsto\ \operatorname{tr}_Lg\,,\qquad g\in G\,.
$$
Since the character is unaffected by conjugation by group elements, we can always restrict it to $\TT$
without any loss of information. We then further consider it as a function  $\chi_L$  on $\mathfrak t$ by
$$
\chi_L :\, x\ \mapsto\ \operatorname{tr}_L \exp({2\pi i\,x})=\sum_\mu\dim L^\mu e^{2\pi i\l\mu,x\r}\,,
$$
where $\mu$ runs over the weights of $L$. In particular, we have the characters $\chi_\lambda=\chi_{L(\lambda)}$, \ $\lambda\in P^+$.

Weyl's character formula is
\begin{equation}
\chi_\lambda(x)=\frac{\sum_{w\in W}(-1)^{l(w)}e^{2\pi i\l w(\lambda+\rho),x\r}}
                  {\sum_{w\in W}(-1)^{l(w)}e^{2\pi i\l w\rho,x\r}}
                  =:\frac{S_{\lambda+\rho}(x)}{S_{\rho}(x)}\,,\qquad
                  \text{for all}\quad x\in\mathfrak t\,.
\end{equation}
Also there is the special product formula for the denominator:
\begin{equation} \label{productFormula}
S_{\rho}(x)=\sum_{w\in W}(-1)^{l(w)}e^{2\pi i\l w\rho,x\r}
           =\prod_{\alpha\in\Delta_+}
           \left(e^{\pi i\l\alpha,x\r}-e^{-\pi i\l\alpha,x\r}\right)\,.
\end{equation}

The functions $S_{\lambda+\rho}(x)$ and $S_{\rho}(x)$ are $W$-skew invariant whereas their quotient is $W$-invariant. The functions 
$S_{\lambda+\rho}(x)= \sum_{w\in W}(-1)^{l(w)}e^{2\pi i\l w(\lambda+\rho),x\r}$ are
called $S$-functions in \cite{KP-S, MP, NPT} due to their similarity in form
to the sine function, which they are in the case of $A_1$. In any case the values of the characters and of the $S$-functions are determined by their values on the fundamental domain $F$, and it
is this fact that becomes the centre of our attention in the sequel.

Note that none of the reflecting hyperplanes \eqref{hyperplane} or \eqref{0-hyperplane}, nor indeed
any reflecting hyperplane of $W_{\mbox{\small{aff}}}$, meets the interior $F^\circ$ of $F$ and so $S_{\rho}(x)$ is never 0 in $F^\circ$. Thus
$|S^2_{\rho}(x)|^2$ is positive on the interior $F^\circ$ of $F$. On the other hand, $S_{\lambda+\rho}(x)$ and $S_{\rho}(x)$ vanish on the boundary of $F$ since $x\in H_j$ implies that $r_jx=x$ while replacing $x$ by $r_jx$ in any $S$-function changes its sign.

\section{$W$-invariant and $W$-skew invariant functions on $\TT$} \label{Winvariance}

\subsection{The algebra of formal exponentials}\ 

Starting with the weight lattice $P$ one may form the algebra $\CC[P]$ of
formal exponentials, which has a $\CC$-basis of symbols $e^\lambda$, $\lambda \in P$,
together with a multiplication defined by bilinear extension of the rule
$$
 e^\lambda \,e^\mu = e^{\lambda +\mu} \,.
 $$
Thus typical elements of $\CC[P]$ are finite complex linear combinations $\sum c_\lambda e^\lambda$.
$\CC[P]$ is an unique factorization domain and its group of invertible elements are
the elements $c e^\lambda$, where $c\ne 0, \lambda \in P$.

 $W$ acts on $P$ and hence as a linear operator and even an automorphism on $\CC[P]$ by $w.e^\lambda = e^{w\lambda}$.
 There is a partial order on $P$ with $\lambda\ge \mu$ if and only if $\mu = \lambda -\beta$
 where $\beta$ is a sum (possibly empty ) of positive roots. This is important because
in the weight systems of irreducible representations of $G$ the highest weight is highest
in this sense. We will use this notion of highest below.
  
There are some advantages to introducing formal exponentials at times since they often clarify the mathematics.
However, we are really interested in their manifestations as functions on $\mathfrak t$ and
on $\TT$ which arise from
\begin{equation}\label{formalToFunction}
x \mapsto e^{2 \pi i x} \mapsto e^{2 \pi i \langle \lambda,x\rangle} \,.
\end{equation}
These mappings are the characters $\phi_\lambda: x \mapsto \exp(2\pi i \lambda.x)$
of the torus, so $\CC[P]$, as functions, is its character ring.
This relates directly back to the previous section where we have defined 
the $G$-characters and $S$-functions, all of which may be viewed as arising
from corresponding elements $\chi_\lambda, S_{\lambda +\rho}$ of $\CC[P]$ . 

At the level of functions the action of $W$ is completely consistent:
$$
w.\phi_\lambda(x) = \phi_\lambda(w^{-1}x) =
\exp(2 \pi i \, \lambda.w^{-1}x) = \exp{2 \pi i \, w\lambda.x} = \phi_{w\lambda}(x)\,.
$$

We are most interested in the subring $\CC[P]^W$ of $W$-invariant elements. 
The simplest forms of $W$-invariant functions are the orbit sums
$\sum_{w\in W} \phi_{w\lambda}$, called $C$-functions in \cite{P,KP-C,MP,NPT}
due to their similarity in form to the cosine function (which they are in type $A_1$).
More relevant here are the characters $\chi_\lambda$ of $G$ (restricted to $\TT$) already introduced in \S\ref{Lie groups}:
\begin{equation}\label{character}
\chi_\lambda = \sum_{\mu \in P} \dim L(\lambda)^\mu e^{2 \pi i \mu}
\end{equation}
where $L(\lambda)$ is the irreducible representation of highest weight $\lambda \in P^+$
and $\dim L(\lambda)^\mu$ is the dimension of its $\mu$-weight space. These are $W$-invariant and they form a basis
for $\CC[P]^W$. 

\begin{prop} \cite{B}{\rm  (Ch.VI)} \quad
$\CC[P]^W$ is a polynomial ring with the fundamental characters $\chi_{\o_1},\dots \chi_{\o_n}$ as the generators
\footnote{The orbit sums over the various dominant weights also form a 
basis for $\CC[P]^W$ and those for $\o_1, \dots, \o_n$ also form a set of generators
for it as a polynomial ring. The relationship between orbit sums and corresponding
characters is a triangular matrix of integers with $1$s down the diagonal \cite{BMP}.}.
\end{prop}

This result really underlies the results of this paper. It says that the fundamental characters of $G$ can be used as new variables by which the algebra of invariant elements of $\CC[P]$
becomes a polynomial ring in these variables. It is in working out the Fourier and functional analysis implied by this statement that the cubature formulae arise.

In the sequel we prefer to reserve the word {\em character} for the characters of $G$
(as opposed to the characters of $\TT$) since they are of fundamental importance to the paper. 

\subsection{Skew invariants elements of $\CC[P]$}\ 

Elements $\xi \in \CC[P]$ for which $w.\xi = (-1)^{l(w)} \xi$ for all $w\in W$ are
called {\em skew-invariants}.  They play a vital role in the paper. The simplest example is $S_\rho$, and it is the foundation for all the skew-invariant elements.

\begin{prop} \cite{B} {\rm(Ch.VI)} \label{skewInvThm}  \quad 
$S_\rho \,\CC[P]^W$ is the set of all $W$-skew-invariant elements of $\CC[P]$.
\end{prop}

Later on, when we use the basic characters $\chi_{\om_j}$ as new variables $X_j$ and have
polynomial functions of the $X_j$, we shall
have need of a Jacobian for the switch of variables from the $X_j$ back to the variables
that parameterize $\TT$. We establish the key result here.

For each $\check \alpha \in \check Q$ there is a unique derivation\footnote{A derivation of
$\CC[P]$ is a linear mapping $D: \CC[P]\longrightarrow \CC[P]$ satisfying 
$D(\phi \xi) = D(\phi) \xi + \phi D(\xi)$
for all $\phi, \xi\in \CC[P]$.}
 $D_{\check\alpha}$
on $\CC[P]$ satisfying  
$$
D_{\check\alpha}(e^\lambda) = \langle\lambda, \check\alpha\rangle e^\lambda \,
$$
for all $\lambda$. $D_{\check\alpha}$ is a linear in $\check\alpha$.

Let $\chi_{\om_1}, \dots, \chi_{\om_n}$ be the basic characters and let
$\check\alpha_1, \dots, \check\alpha_n$ be the standard basis of 
$\check Q$ dual to $\{\om_1, \dots, \om_n\}$. Let $J$
be the matrix with entries $J_{jk} = D_{\check\alpha_j}\chi_{\om_k}$.

\begin{prop} {\rm (Steinberg\footnote{This result is an Exercise to Ch. VI in \cite{B} and
 attributed to R. Steinberg there.})} \label{Steinberg}
$$
\det(J) = S_\rho\,.
$$
\end{prop}
\begin{proof} All the exponentials in $\chi_{\om_k}$ are of the form $e^{\om_k - \beta_k}$
where $\beta_k$ is a sum of positive roots, and the highest term is $e^{\om_j}$. 
Now $D_{\check\alpha_j}\chi_{\om_k}$ is a sum of exponentials of the same form, but
since $\langle\om_j, \check\alpha_k\rangle =\delta_{jk}$, the highest terms only survive
along the diagonal of $J$. Thus when we compute the determinant we obtain a sum
of signed products $e^{\om_1 -\beta_1} \dots e^{\om_n -\beta_n}$ and only the term 
from the diagonal can contribute an exponential of the form $e^{\om_1} \dots e^{\om_n }$, 
and its coefficient is $1$.
Thus $\det(J)$ is an element of $\CC[P]$ whose highest term is $\om_1 +\cdots + \om_n = \rho$
and this occurs with coefficient equal to $1$.

We shall prove that $\det(J)$ is $W$-skew invariant. Then by Prop.~\ref{skewInvThm} it is a 
multiple of $S_\rho$. Because all the weights in the expansion of $\det(J)$ are less than or equal
$\rho$, this multiple can only be a scalar. Since the leading coefficient is $1$ in both
cases, $\det(J) = S_\rho$, as we wish to prove. 

A simple computation shows that $wD_{\check\alpha} w^{-1} = D_{w\check\alpha}$ for
all $\check\alpha$. Fix any $l =1, \dots, n$ and let $r_l= r_{\check\alpha_l}$. 
We have $r_l\check\alpha_j = \check\alpha_j - A_{lj} \check\alpha_l$.
Then,
\begin{eqnarray*}
r_l \det(J) &=& \det(r_l(J)) = \det((r_l (D_{\check\alpha_j}\chi_{\om_k})))\\
&=& \det((r_l D_{\check\alpha_j}r_l r_l\chi_{\om_k})) = \det((D_{r_l\check\alpha_j} \chi_{\om_k}))\\
&=& \det((D_{\check\alpha_j - A_{lj} \check\alpha_l}(\chi_{\om_k})))
= \det((D_{\check\alpha_j} - A_{lj} D_{\check\alpha_l})(\chi_{\om_k})))\\
&=& \det((D_{\check\alpha_j}(\chi_{\om_k}) - A_{lj} D_{\check\alpha_l}(\chi_{\om_k})))\,,
\end{eqnarray*}
where in the second line we used the $W$-invariance of the characters. 
The operation has resulted in altering all rows except the $l$th row by a multiple of
the $l$th row (which does not alter the determinant) and replacing the $l$th row by
its negative, which changes the sign of the determinant. Thus $r_l \det(J) = - \det(J)$,
which gives the desired skew-symmetry.
\end{proof}

\subsection{An inner product on $\CC[P]^W$ } \ 

The natural inner product on $\CC[P]$ is $\langle\cdot\,,\cdot\rangle_\TT$ defined by
$$
\langle f, g\rangle_\TT = \int_\TT f \overline g \, d\theta_\TT\,,
$$
where the $\theta_\TT$ is the normalized Haar measure of the torus. Relative to this
the functions $\phi_\lambda$ form an orthonormal basis. 
Using this inner product we can
complete $\CC[P]$ in the corresponding $L^2$-norm to the Hilbert space $L^2(\TT,\theta_\TT)$
with the normalized $\phi_\lambda$ forming an orthonormal basis, in the sense of Hilbert spaces. 
Of course we can look at the closure of $\CC[P]^W$ in $L^2(\TT,\theta_\TT)$, which
is in fact the subspace of $W$-invariant elements  $L^2(\TT,\theta_\TT)^W$ of
$L^2(\TT,\theta_\TT)$.

However, the inner product $\langle \cdot\,,\cdot\rangle_\TT$ is not ideal for this subspace, and rather we would like to find one with respect to which the characters $\chi_\lambda$ form an orthonormal base. 

We note that for any $f\in L^2(\TT,\theta_\TT)^W$, $f S_\rho \in L^2(\TT,\theta_\TT)$,
and it is skew-invariant with respect to $W$. Form its Fourier expansion
$$ 
fS_\rho = \sum_\mu \langle f S_\rho, \phi_\mu \rangle_\TT \, \phi_\mu \,,
$$
equality being in the $L^2$ sense. 
The summands can be gathered together into $W$-orbits, and on each orbit the coefficients
$\langle f S_\rho, \phi_\mu \rangle_\TT$ are equal in absolute value and alternate
in sign according to the parity of the Weyl group elements. The only orbits
that do not vanish are those contain an weight $\mu \in P^{++}$, and we get
$$ 
fS_\rho = \sum_{\lambda \in P^+} \langle f S_\rho, \phi_{\lambda+\rho} \rangle_\TT \, S_{\lambda + \rho} \,.
$$
Dividing out the function $S_\rho$, which is valid as long as the functions are restricted to
the interior of the fundamental chamber $F^\circ$, we obtain
\begin{equation} \label{FEoff}
 f = \sum_{\lambda \in P^+} \langle f\,S_\rho, \phi_{\lambda+\rho}\rangle_\TT \,
\chi_{\lambda} \, .
\end{equation}
Now, using the $W$-invariance of $\theta_\TT$ and the skew-invariance of
$fS_\rho$, we have
\begin{eqnarray}\label{rewriteInnerProduct}
 \langle fS_\rho, \phi_{\lambda+\rho}\rangle_\TT &=& 
 \int_\TT f S_\rho \,\overline{\phi_{\lambda+\rho}} \, d\theta_\TT 
 = \frac{1}{|W|} \int_\TT \sum_{w\in W}(-1)^w f S_\rho \,
 \overline{\phi_{w(\lambda+\rho)}} \, d\theta_\TT \nonumber\\
 & =&   \frac{1}{|W|} \int_\TT f \,S_\rho \overline{S_{\lambda+\rho}}\, d\theta_\TT
= \frac{1}{|W|} \int_\TT f \,
 \overline{\chi_{\lambda}} \,S_\rho \overline{S_\rho}\, d\theta_\TT
 = \int_{F^\circ} f \,
 \overline{\chi_{\lambda}} \,S_\rho \overline{S_\rho}\, d\theta_\TT\,.
 \end{eqnarray}
 
 Thus \eqref{FEoff} and \eqref{rewriteInnerProduct} show that $f\in L^2(\TT,\theta_\TT)^W$ has a Fourier expansion in terms
 of the characters $\chi_{\lambda} $, $\lambda \in P^+$, with 
 coefficients given by a new inner product defined on $L^2(\TT,\theta_\TT)^W$
 by 
 \begin{equation}\label{newInnerProduct}
 (f, g) \mapsto  \int_{F^\circ} f \, \overline{g} \,S_\rho \overline{S_\rho}\, d\theta_\TT\,.
 \end{equation}
 We can then rewrite \eqref{FEoff} as
 \begin{equation} \label{newFEoff}
 f = \sum_{\lambda \in P^+} ( f ,\chi_\lambda) \,\chi_\lambda \,.
 \end{equation}

 We shall use these results in \S\ref{Approximations}.

%
\section{Elements of Finite Order}\label{EFOs}
\medskip

Elements of finite order \cite{MP84} (EFOs) in $G$ are used to create the interpolation points for
the discrete Fourier analysis  \cite{MP87} and cubature formulae to follow. 

We have seen that every element of $G$ is conjugate to one of the form $g=\exp(2 \pi i x)$
where $x\in F$. The element $g$ is called {\em regular} if its centralizer is of dimension $n$,
the rank of $G$. Since each element lies in a torus this is the smallest possible dimension for a centralizer.
Regularity is a property of the entire conjugacy class of an element
and for $x\in F$ it is equivalent to saying that $x \in F^\circ$.

The condition that $g$ has finite order dividing $M$ is the equivalent to the condition
that $\exp(2 \pi i x)^M = \exp(2 \pi i Mx)$ acts trivially on every irreducible representation,
and for this all we need is that it acts trivially on every weight space $L(\lambda)^\mu$.
In turn this requires precisely that $\l \mu, Mx\r \in \ZZ$ for all weights of $P$,
and finally it is equivalent to $M x \in \check Q$, since $\check Q$ is the $\ZZ$-dual
of $P$. 

In fact what we are going to need here is not that $x \in  \frac{1}{M} \check Q$ but rather that 
$$ 
x \in  \frac{1}{M} \check P\, ,
$$
a statement that is equivalent to saying that $Ad(g)^M =1$,
i.e. $g^M$ acts trivially in the adjoint representation.  In this case we say that $g$ has $Ad$-{\em order}
or {\em adjoint order} $M$, even though the actual adjoint order, which we shall call
the {\em strict} adjoint order, namely the {\em least} $N$ for which $Ad(g)^N =1$ 
may be some proper divisor $N$ of $M$. We also say that $x\in\mathfrak t$ is an element of
adjoint order $M$ if $\exp(2\pi i x)$ is of adjoint order $M$.

Given the definition above, the conjugacy classes of elements of adjoint order $M$ are represented by the points $x$ of the form
\begin{equation}\label{point}
x=  \frac1M\sum_{j=1}^ns_j\check\omega_j
\end{equation}
where
\begin{equation} \label{AdOrderMconditions}
s_j \in \ZZ^{\ge 0} \quad \mbox{for all } j \quad \mbox{and} \quad
\l-\alpha_0,x\r = \frac1M \sum_{j=1}^n m_j s_j \le 1\,.
\end{equation}
The regular conjugacy classes of adjoint order $M$ are 
represented by (\ref{AdOrderMconditions}) where the inequalities
are made strict. We write $F_M$ (resp. $F^\circ_M$) for the elements
of $F$ (resp. $F^\circ$) of $Ad$-order $M$.

Using $m_0=1$ defined in \S\ref{Lie groups} we can define $s_0 \in \ZZ^{\ge0}$
so that 
\begin{equation}\label{KacCoordCond}
 \sum_{j=0}^n m_j s_j =M \,.
 \end{equation}
Listing all the elements of $F_M$ (resp. $F^\circ_M$) is then just a question of finding
all non-negative (resp. positive) integer solutions $[s_0, s_1, \dots, s_n]$ to (\ref{KacCoordCond}). 
We call $[s_0, s_1, \dots, s_n]$ the {\em Kac coordinates} of $x$. 

We will be particularly interested in the set $F_{M+h}$ of elements $x\in F$ of $Ad$-order $M+h$ for some non-negative integer $M$:
\begin{equation} \label{adOrderM+h}
x=\tfrac{1}{M+h} (s_1\check\o_1+s_2 \check\o_2+\cdots+s_n\check\o_n)
   \in\tfrac1{M+h} \check P\,,\qquad\text{where}\quad \sum_{j=1}^nm_js_j\leq M+h
\end{equation}
where $s_1\geq0,\dots,s_n\geq0$ are integers such that $\sum_{j=1}^nm_js_j\leq M+h$. Alternatively we have the Kac coordinates $[s_0,s_1,\dots,s_n]$. 

\bigskip

Each of the following three conditions assures that $x$ of \eqref{adOrderM+h} is in 
$F^\circ_{M+h}$:
\begin{equation}         
\begin{alignedat}{2}\label{intCondition}
&s_j>0 \,,                & \qquad &j=0,1,\dots,n\,;\\
&\sum_{j=0}^n m_jt_j=M\,, &\qquad  &t_j:=s_j-1\geq0\,; \\
&\sum_{j=1}^nm_jt_j\leq M &\qquad  &\mbox{(so $t_0$ completes the sum to $M$)}\,.
\end{alignedat}
 \end{equation}
 \medskip

When $M=0$, it contains only the element given by $s_0=s_1=\cdots=s_n=1$. For $M \ge 0$, it clearly contains $|F^\circ_{M+h} |=|F_{M}|$ points.
Formulas for the cardinality of $|F_{M+h}|$ have been worked out for all $M$ and for all simple $G$ in \cite{HP}.

\section{Points of $F^\circ_{M+h}$ as zeros of $S$-functions}\label{start}

We are now at a point where we begin the main development of the paper. We fix, once
and for all a non-negative integer $M$. The first
step is to show that the points of $F^\circ_{M+h}$ are common zeros of a certain
set of $S$-functions. These points are the interpolation points for the cubature
formulae to follow.

\medskip
Consider a dominant weight $\lambda= \lambda_1\o_1 + \dots + \lambda_n \o_n$ \ .
We want to find points $x\in F^\circ_{M+h} $ at which the $S$-function 
$S_{\lambda+\rho}(x)$ vanishes:
\begin{equation*}
S_{\lambda+\rho}(x)=\sum_{w\in W}(-1)^{l(w)}
                   e^{2\pi i\l w(\lambda+\rho),x\r} =\sum_{w\in W}(-1)^{l(w)}
                   e^{2\pi i\l\lambda +\rho,w^{-1}x\r}=0\,.
\end{equation*}

One way to make this happen is to have
\begin{gather}\label{zerosCond}
\l\lambda+\rho,x\r - \l\lambda+\rho,rx\r\in\ZZ\,,\;\mbox{for all} \; x\in\tfrac1{M+h}\check P\, ,
\end{gather}
where $r$ is the reflection in the highest coroot, for if this is the case then
the sum collapses in pairs adding up to zero. Via $\l w\alpha,\check \beta\r =
\l \alpha,w^{-1}\check \beta\r$ we obtain that $r$ appears as the reflection in
some root $\gamma$ on the root/weight side of the picture:
$r x = x - \l \gamma,x\r\check\alpha_0$. 
The condition \eqref{zerosCond} is equivalent to
$$
\l\gamma,x\r\l\lambda+\rho,\check\a_0\r\in\ZZ,\;\mbox{for all} \; x\in\tfrac1{M+h}\check P \, .
$$
Note that $\l\gamma,x\r\in\tfrac1{M+h}\ZZ$ \  since \ $\l\gamma,\check P\r\subset\ZZ$, so we only need  
\ $\l\lambda+\rho,\check\a_0\r\in (M+h)\ZZ$. 

The simplest case is to look for \ 
$\l\lambda+\rho,-\check\a_0\r=M+h$, that is,
\begin{equation*}
\sum_{j=1}^n(\lambda_j+1)\l\o_j,-\check\a_0\r
  =\sum_{j=1}^n(\lambda_j+1)\check m_j = M +h\,,
\end{equation*}
or equivalently,
\begin{equation}\label{EFOsolutions}
\sum_{j=1}^n \lambda_j \check m_j = M+1 \,.
\end{equation}

\medskip

All solutions $\lambda=(\lambda_1,\dots,\lambda_n)$ to \eqref{EFOsolutions}, where the $\lambda_j\in\ZZ^{\geq0}$,  lead to $S$-functions $S_{\lambda+\rho}$ that are zero at all EFOs of $Ad$-order $M+h$ in the interior of the fundamental domain.

\section{Introducing the polynomial functions}

Following \cite{NPT}, assign variables $X_j$ to the characters for weights $\o_j$. Thus we have polynomial variables
\begin{equation}\label{variables}
X_1,\ X_2,\dots,X_n\,,\quad\text{where}\quad
X_j:=\chi_{\omega_j}(x)\,,\qquad x\in F^\circ\,.
\end{equation}

With these we can introduce the domain 
\begin{equation}\label{domain}
\Omega :=\{(X_1(x), \dots, X_n(x)) : x\in F^\circ\} \subset \CC^n\,.
\end{equation}
We shall soon see that this is actually an open subset of a real $n$-dimensional
space and eventually it will be the natural domain of the functions real-valued
functions of the variables $X_1, \dots, X_n$ that we wish
to study.

Define the $m${\em-degree}\footnote{This might be more properly called
the $\check m$-degree, but this seems a bit cumbersome.} of the variables by assigning degree $\check m_j$ to $X_j$.

\begin{table}[h]
\footnotesize
\addtolength{\tabcolsep}{-3pt}

\begin{center}
\begin{tabular}{|c||c|c|c|c||c|c|}
\hline
variable &$X_1$  
          &$X_2$   
          &$\cdots$
          &$X_n$
          &$S_\rho$
          &$K$  \\
\hline 
$m$-degree   &$\check m_1$ 
             &$\check m_2$
             &$\cdots$
             &$\check m_n$
             &$h-1$
             &$2h-2$  \\
\hline
\end{tabular}

\bigskip

\caption{The $m$-degrees of the polynomial variables $X_1,\dots,X_n$ and of the functions $S_\rho$ and $K$. The Coxeter number $h=1+\check m_1+\cdots+\check m_n$.} \label{m-deg}
\end{center}
\end{table}
\noindent
Then the monomials $X_1^{\lambda_1}\dots X_n^{\lambda_n}$ \ of $m$-degree \ $\leq M$ \ are those satisfying 
\begin{equation}\label{m-degreeInequal}
\lambda_1\check m_1+\cdots+\lambda_n \check m_n\leq M
\end{equation}
where \ $\lambda_1\geq0,\dots,\lambda_n\geq0$. Although the marks and co-marks
are not necessarily identical,  Fig.~\ref{diagrams}, they are at worst simply permutations of each other. Thus 
\eqref{m-degreeInequal} has the same number of solutions $(\lambda_1,\dots,\lambda_n)$ as we saw before, namely $|F^\circ_{M+h}|=|F_{M}|$. 
The constant polynomials are those of $m$-degree $0$.

In keeping with this notation, we will say also that  $\lambda=(\lambda_1,\dots,\lambda_n)$ has $m$-degree equal to 
\begin{equation}\label{m-degree}
\lambda_1\check m_1+\cdots+\lambda_n\check m_n = \l \lambda, -\check\a_0\r \,.
\end{equation}

\begin{theorem}\label{0points}\ 

The number of monomials $\CC[X_1,\dots,X_n]$ of $m$-degree $\leq M$ is equal to the number of regular EFOs of $Ad$-order $M+h$ in the fundamental chamber. Each of the  regular EFOs of $Ad$-order $M+h$ in the fundamental chamber is a common zero
of all the $S$-functions $S_{\lambda+\rho}$ and all the character functions
$\chi_\lambda$ for which $\lambda$ has $m$-degree equal to $M+1$.  \qed
\end{theorem} 

The trick that we have
used above of using the internal reflective anti-symmetry to construct common zeros  is taken from \cite{LX}. It is remarkable that in the case of type $A_n$ root lattices it actually finds {\em all} the common zeros. The proof of this makes essential use of the fact that the new variables
$X_1, \dots, X_n$ are all of degree $1$, something that is true only for type $A_n$. In fact our example
of type $G_2$ in \S\ref{G2} indicates that this result does not hold in general.
\smallskip

The `smallest' $S$-function is the one defined by the strictly dominant weight 
of lowest $m$-degree, namely $\rho$ of \eqref{rho} with $m$-degree $h-1$.
Writing $S_\rho$ in its well-known form \eqref{productFormula}, 
we note that $\overline{S_\rho}
= (-1)^{|\Delta_+|} S_\rho $. Thus 
$$
|S_\rho|^2 = S_\rho \overline S_\rho =  (-1)^{|\Delta_+|} S_\rho^2 \, ,
$$
 and we note that this function is positive on all of $F^\circ$ and vanishes on its boundary. $S_\rho \overline S_\rho$ is 
a $W$-invariant function and so is expressible as a polynomial in the basic
characters $\chi_{\o_j} = X_j$. We then have the corresponding strictly positive function 
$K$ on $\Omega$:
\begin{equation} \label{theK}
K(X_1, \dots, X_n) = S_\rho \overline S_\rho(x),\quad x \in F^\circ \,.
\end{equation} 

We note here that if $\mu \in P^{++}$ then 
$$
\overline S_\mu =
\sum_{w\in W} (-1)^{l(w)} e^{- w\mu} = \sum_{w\in W} (-1)^{l(w)} (-1)^{l(w_{opp})} e^{-w w_{opp} \mu}
= (-1)^{|\Delta_+|}\,S_{-w_{opp} \mu} \,,
$$
where $w_{opp}$ is the opposite involution in $W$ (since $w_{opp}$ is a product of $ |\Delta_+|$
reflections). 

The opposite involution interchanges the positive roots (resp. positive coroots) with the negative ones, and, since $\mu$ is dominant so is $-w_{opp}\, \mu$. Of course $w_{opp}$ is not always simply the negation operator, see Tab.~\ref{permute}. Still, it does simply change the sign of the highest positive root (resp. highest coroot). Thus we have the important little equation
\begin{equation}\label{wopp}
m\mbox{-deg}\,(\mu) = \l \mu,-\check\a_0\r  = \l w_{opp}\, \mu, -w_{opp} \,\check\a_0 \r 
= \l w_{opp}\, \mu, \check\a_0 \r = m\mbox{-deg}\,(-w_{opp}\, \mu) \, .
\end{equation}

This is useful because it means that $(-1)^{|\Delta_+|}\,\overline{S_\mu}$ and $\overline{\chi_\mu}$ are just another $S$-function and another group character respectively,
and the highest weight involved in each case has the same $m$-degree as before conjugation. In 
particular conjugation of the characters $\chi_{\om_j}$ can at worst permute some of them, say
$\chi_{\om_j} \mapsto \chi_{\sigma(\om_j)}$
by some permutation $\sigma$ of order $2$ of the indices $\{1, \dots, n\}$. 
Table~\ref{permute} shows what happens in the cases when $\sigma$ is not just the identity
permutation. 

If we let 
\begin{equation}
\mathfrak R =\{ z= (z_1, \dots, z_n): \overline{z} = (z_{\sigma(1)}, \dots,
z_{\sigma(n)})\}\,,
\end{equation}
then $\mathfrak R$ is a real space of dimension $n$ 
and $\Omega \subset \mathfrak R$ is an $n$-dimensional subdomain, as follows
from the non-vanishing of the Jacobian $S_\rho$ on $F^\circ$ \eqref{jac}. This is the space
on which we shall think of $X_1,\dots, X_n$ as real variables.

\begin{table}[h]
\footnotesize
\addtolength{\tabcolsep}{-3pt}

\begin{center}
$A_n$\ \ $(n>1)$\qquad
\begin{tabular}{|c||c|c|c|c|c|c|}
\hline
$X$       &$X_1$  
          &$X_2$   
          &$\cdots$
          &$X_{n-1}$
          &$X_n$
           \\
\hline 
\rule{0pt}{10pt}
$\overline X$   &$X_n$ 
                &$X_{n-1}$
                &$\cdots$
                &$X_2$
                &$X_1$
          \\
\hline
\end{tabular}

\bigskip

$D_{2n+1}$\ \ $(n>1)$\qquad\ \quad
\begin{tabular}{|c||c|c|c|c|c|c|}
\hline
$X$       &$X_1$  
          &$X_2$   
          &$\cdots$
          &$X_{2n-1}$
          &$X_{2n}$
          &$X_{2n+1}$
           \\
\hline 
\rule{0pt}{10pt}
$\overline X$   &$X_1$ 
                &$X_2$
                &$\cdots$
          &$X_{2n-1}$
          &$X_{2n+1}$
          &$X_{2n}$
          \\
\hline
\end{tabular}

\bigskip

$E_6$\qquad\qquad\qquad
\begin{tabular}{|c||c|c|c|c|c|c|}
\hline
$X$       &$X_1$  
          &$X_2$   
          &$X_3$
          &$X_4$
          &$X_5$
          &$X_6$
           \\
\hline 
\rule{0pt}{10pt}
$\overline X$   &$X_5$ 
                &$X_4$
                &$X_3$
          &$X_2$
          &$X_1$
          &$X_6$
          \\
\hline
\end{tabular}

\bigskip
\caption{Correspondence of the variables $X_j$ to $\overline X_j$, $(j=1,\dots,n)$ produced by the action of $-w_{opp}$. In all other cases $\overline{X_j}=X_j$.} \label{permute}
\end{center}
\end{table}
\noindent

\begin{remark}\label{conjugation!}
As we have seen, conjugation actually permutes some of the basic variables $X_j$. 
We shall use the overline symbol to indicate this form of conjugation. Thus one should understand
the conjugation symbol as having this dual meaning of actual complex conjugation when
the $X_j$ are treated as functions on $\TT$ and as the permutation $\sigma$ when treated
as the coordinate variables of $\mathfrak R$.
Thus we shall write $\overline{c\, X_{j_1}\dots X_{j_r}}$ (where $c \in \CC$) to mean
$\overline c \, \overline{X_{j_1}} \dots \overline{X_{j_r}}$, understanding that 
the $\overline X_j$ has this dual meaning. For a polynomial $g(X_,\dots, X_n) =
\sum c_{j_1,\dots,j_n} X_1^{j_1}\dots X_n^{j_n}$, $\overline{g(X_1,\dots, X_n)} :=
\sum \overline{c_{j_1,\dots,j_n}} \overline{X_1}^{j_1}\dots \overline{X_n}^{j_n}=
\overline g( \overline{X_1}^{j_1}\dots \overline{X_n}^{j_n})$.
\end{remark}

Notice that since
$$
X_j \overline{X_k} K\longleftrightarrow \chi_{\omega_j} \overline{\chi_{\omega_k }}S_\rho \overline{S_\rho} = S_{\omega_j+\rho} \overline{S_{\omega_k+\rho}}
$$
and
$$\langle\omega_j +\rho,-\check\alpha_0\rangle =\check m_j+h-1 \quad \mbox{and} \quad 
\langle-w_{opp}(\omega_k +\rho),-\check\alpha_0\rangle =\check m_k+h-1 \, ,$$
we understand that $K$ has $m$-degree equal to $2h-2$.

\section{The integration formula}\label{The Integration Formulae}

We wish to study weighted integrals of the form
$$
\int_\Omega  f\, \overline g \,K^{1/2} dX_1\dots dX_n\,,
$$
where $f,g$ are functions of the variables $X_1,\dots, X_n$ defined on $\Omega$.
These are related back to $\mathfrak t$ (more specifically to $F^\circ$) and the torus $\TT$ via the defining equations \eqref{variables}.

\subsection{The key integration formula}\ 

Natural variables for $\mathfrak t$  are  $x = (x_1, \dots, x_n) =  \sum_{j=1}^n x_j \check\alpha_j$ where
the $x_j$ run over $[0,1) \times \dots \times [0,1)$. 

 The derivation $D_{\check\alpha_j}$ on $\CC[P]$, 
$ D_{\check\alpha_j}e^\lambda = \langle\lambda, \check\alpha_j\rangle e^\lambda$ is, when $\CC[P]$ is treated as an algebra
 of functions on $\TT$, the mapping
 \begin{eqnarray}D_{\check\alpha_j} e^{\langle\lambda, 2 \pi i x\rangle} &=& 
 D_{\check\alpha_j} e^{\langle\lambda, 2 \pi i \sum x_k\check\alpha_k\rangle} \nonumber\\
 & = &
 \langle\lambda, \check\alpha_j\rangle e^{\langle\lambda, 2 \pi i x\rangle}
 = \frac{1}{2\pi i}\frac{d}{dx_j}e^{\langle\lambda, 2 \pi i x\rangle} \,.
 \end{eqnarray}
 
 Using Prop.~\ref{Steinberg} we then see that the Jacobian of the transformation of the variables $x$ to 
 variables $X$ is 
 \begin{equation}\label{jac}
 |(2 \pi i)^n S_\rho(x) | = (2 \pi )^n | S_\rho(x) | \,.
 \end{equation}
 Thus from the definition of $K$ we have 
\begin{gather*}
\int_\Omega  f\, \overline g\,K^{1/2} dX_1\dots dX_n = \int_\Omega  f(X_1,\dots, X_n) \, \overline {g(X_1,\dots, X_n)}\,K^{1/2}(X_1,\dots,X_n)\,dX_1\dots dX_n \\
=(2\pi)^n\int_{F^\circ}f(\chi_{\o_1}(x),\dots,\chi_{\o_n}(x))\,
{\overline g(\overline{\chi_{\o_1}}(x),\dots,\overline{\chi_{\o_n}}(x))}\,
S_\rho(x)\,\overline{S_\rho}(x)\,dx_1\dots dx_n
\end{gather*}
for all functions $f,g$ are
in the variables $X_1, \dots, X_n$ on $\Omega$.

\begin{theorem} \label{integralFormula}\ 

Let $M$ be a positive integer. Then for all polynomials $f,g \in \CC[X_1,\dots, X_n]$
with $m\mbox{-deg}\,(f) \le M+1$ and $m\mbox{-deg}\,(g) \le M$ we have
\begin{equation}\label{intForm}
\begin{aligned}
&\int_\Omega f\,\overline g\,K^{1/2}\,dX_1\dots dX_n \\
&=(2\pi)^n\int_{F^\circ}f(\chi_{\o_1}(x),\dots,\chi_{\o_n}(x))\,
\overline {g(\chi_{\o_1}(x),\dots,\chi_{\o_n}(x))}\,
S_\rho(x)\,\overline{S_\rho}(x)\,dx_1\dots dx_n  \\
&=\frac{1}{c_G} \left(\frac{2 \pi }{M+h} \right)^n\sum_{x\in {F^\circ_{M+h}}}
   f(\chi_{\omega_1}(x),\dots,\chi_{\omega_n}(x))
   \overline {g(\chi_{\o_1}(x),\dots, \chi_{\o_n}(x))} \,
   S_\rho(x)\overline{S_\rho}(x) \, .
\end{aligned}
\end{equation}
\end{theorem}

This theorem is proved in \S\ref{intFormProof}.

\subsection{The cubature formula}\ 

For any function $f$ defined on $\Omega$, let $\widetilde f$ be defined by 
$\widetilde f (x) = f(\chi_{\omega_1}(x),\dots,\chi_{\omega_n}(x))$.

\begin{coro}\label{cubature}\ 
 
Let $M$ be a non-negative integer. Then for all polynomials $f \in \CC[X_1,\dots, X_n]$
with $m\mbox{-deg}\,(f) \le 2M+1 $ we have
\begin{equation}\label{cf}
\int_\Omega f\,K^{1/2}\,dX_1\dots dX_n = 
\frac{1}{c_G} \left(\frac{2 \pi }{M+h} \right)^n\sum_{x\in {F^\circ_{M+h}}}
   \widetilde f(x)\,\widetilde K(x) \, .
\end{equation}
\end{coro}
Equation \eqref{cf} is the {\bf cubature formula}. The points $x \in {F^\circ_{M+h}} =\left(\tfrac{1}{M+h}\check P\right)\cap F^\circ$ are the interpolation points. These points are common zeros of the character functions of the $m$-deqree $M+1$. The coefficients of the interpolation are the values of 
$\widetilde K = |S_\rho|^2$ at the interpolation points. The Corollary is a direct consequence of
Theorem~\ref{integralFormula} since every polynomial of $m$-degree less than or equal to 
$2M+1$ can be written as a linear combination of polynomials of the form $f\overline g$ appearing in \eqref{intForm}.

We first prove a key lemma.

\subsection{Separation lemma}\

\begin{lemma} \label{keyLemma}\ 

If $\phi = \sum_{j=1}^n \phi_j \o_j \in P^+$ and $\phi \ne 0$, and if
 $m\mbox{-deg}\,(\phi) < 2(M+h)$, then  $\phi \notin (M+h)\,Q$.
\end{lemma}
\begin{proof} Suppose by way of contradiction that $\phi \in (M+h)\,Q$. We have
$$
0 < \l\phi,-\check\a_0\r < 2(M+h)
$$
from our assumption on the $m$-degree of $\phi$ and since $\l\phi,-\check\a_0\r = \sum \phi_j\,\check m_j >0$.
However, since $\l Q,-\check\a_0\r \subset \ZZ$, our assumption on $\phi$ forces 
$\sum \phi_j\,\check m_j = \l\phi,-\check\a_0\r \in (M+h)\,\ZZ$ and hence $\sum \phi_j\,\check m_j =M+h$.

In the same way, applying $\phi$ to each simple coroot $\check\a_j$ in turn, we obtain
$\phi_j=\l \phi, \check\a_j\r = (M+h)\,a_j$ for some $a_j \in \ZZ^{\ge 0}$, $j=1, \dots, n$.
Thus $M+h = \sum \phi_j\,\check m_j = (M+h)\sum_{j=1}^n a_j\, \check m_j$. This implies that
exactly one $a_j \ne0$ and for this $j$, $a_j=1$ and $\check m_j=1$. Thus
$\phi = (M+h)\o_j \in (M+h)\,Q$, i.e. $\om_j \in Q$.

In fact this can't happen. One way to see this is to use a well-known fact that the non-trivial
elements of the centre of the simply-connected simple  Lie group $G$ are 
given by the elements $\check\omega_j$ over the $j >0$ which correspond the places
where $m_j=1$ (these are certain vertices, different from the vertex $0$, of the fundamental 
chamber. Of course these elements $\check\omega_j \notin \check Q$ for such
an element would be a representative of the identity element. Our case here is the same situation except it is for the simply-connected simple  Lie group $\check G$ based on the dual
root system: we have 
a $j$ for which $\check m_j=1$, and hence $\omega_j  \notin Q$.

Thus $\phi = (M+h)\o_j$ cannot lie in $(M+h)\,Q$.
\end{proof}

\subsection{Weyl integral formula and its consequences}\ 

We recall here the Weyl integral formula \cite{BS}, 
$$
\int_G \mathcal F\, d\theta_G =
\frac{1}{|W|} \int_\TT \mathcal F\, |S_\rho|^2 \,d\theta_\TT
$$
for all class functions $\mathcal F$ (functions which are invariant on conjugacy classes in $G$). Here the measures are normalized Haar measure on $G$ and $\TT$ respectively and the function $\mathcal F$ is simply being restricted to the maximal torus $\TT$ in the second integral. In particular characters are class functions, and so for the  irreducible characters $\chi_\lambda, \chi_\mu$ of the irreducible representations of $G$ of highest weights $\lambda, \mu$ respectively, we have from the standard orthogonality relations \cite{BS}:
\begin{equation} \label{wif}
\delta_{\lambda,\mu} = \int_G \chi_\lambda \overline{\chi_\mu} d\,\theta_G =
\frac{1}{|W|} \int_\TT \chi_\lambda \overline{\chi_\mu} \, |S_\rho|^2 \,d\theta_\TT\,.
\end{equation}
Here we are using the usual Kronecker delta. 

We do not need all of \eqref{wif}, only the equality of the left and
right hand sides; and that fact is not hard to see. We have
$$
\chi_\lambda S_\rho \overline{\chi_\mu S_\rho} = S_{\lambda +\rho} \overline{S_{\mu+\rho}}
= \sum_{w\in W} \sum_{v\in W} (-1)^{l(w)}(-1)^{l(v)}\, e^{2\pi i \l w(\lambda+\rho),x\r)} 
e^{-2\pi i  \l v(\mu+\rho),x\r)}\,
$$
by Weyl's character formula\footnote{In \cite{BS} Weyl's character formula is derived
from the integral formula, but algebraists usually use an algebraic proof of the
character formula.}.
The integral over $\TT$ of $ e^{2\pi i \l w(\lambda+\rho),x\r)-2\pi i  \l v(\mu+\rho),x\r)}$
is $0$ unless $w(\lambda+\rho)- v(\mu+\rho) =0$, in which case it integrates to $1$. Since $\lambda+\rho$ and
$\mu+\rho$ are strictly dominant, this happens only if $\lambda = \mu$ and
$w=v$. If indeed $\lambda = \mu$ then there are exactly $|W|$ times when
$w(\lambda+\rho)- v(\mu+\rho) =0$, and we see that the right hand side of 
\eqref{wif} is $\delta_{\lambda,\mu}$.


We also wish to recall a result from discrete Fourier analysis \cite{MP}. There it is 
proved that, using the notation established above, 
\begin{equation}\label{dfa}
\int_F S_{\lambda+\rho} \overline{S_{\mu+\rho}}\, d\theta_{\TT}
= \frac1{c_G(M+h)^n}\sum_{x\in {F^\circ_{M+h}}} S_{\lambda+\rho}(x) \overline{S_{\mu+\rho}(x)}
=\delta_{\lambda, \mu}\,,
\end{equation}
as long as when $\lambda \ne \mu$, the points of $F^\circ_{M+h}$, can {\em separate} all the weights appearing
in the $W$-orbit of $\lambda+\rho$ from all those appearing in the $W$-orbit
of $\mu+\rho$.  Explicitly, separation means that it never happens that
$w(\lambda+\rho)-v(\mu+\rho)$ takes integer values on all the points of 
$\tfrac{1}{M+h} \check P$, or equivalently, it never happens that
 $w(\lambda+\rho)-v(\mu+\rho) \in (M+h)\,Q$, except when $\lambda = \mu$.

This proof of the last equality in \eqref{dfa} is actually a straightforward thing to see. First note that 
$S_{\lambda+\rho}\, \overline{S_{\mu+\rho}}$ is $W$-invariant, and hence its integral
over all of $\TT$ is $|W|$ times its integral over $F$. There is no need to worry about the
boundary of $F$ which has measure $0$ and in any case the function takes the value
$0$ on all of its boundary. In the same way, the sum over $F^\circ_{M+h}$ can be extended by
the operation of the Weyl group to obtain a full set of representatives of $\tfrac{1}{M+h} \check P/ \check Q$, noting again that since the function is $0$ on the boundaries of the chambers, adding in the extra boundary elements that may appear in $\tfrac{1}{M+h} \check P/ \check Q$ makes no difference to the sum. This again increases the value of the sum by $|W|$ but has the benefit of turning the sum into a sum over a group.  Then usual considerations of sums of exponentials over a group give the desired orthogonality  as long as the separation condition is satisfied. The order of $\tfrac{1}{M+h} \check P/ \check Q$ is $c_G\,(M+h)^n$, which explains the factor outside the sum part of the formula.
\smallskip

\subsection{Demonstration of Theorem~\ref{integralFormula}}\label{intFormProof} \

With these facts, we now prove the integral formula \eqref{intForm}:
\begin{proof} 
Because of linearity, it is enough to prove \eqref{intForm} when $f,g$ are monomials
of the form
\begin{gather*}
\chi_{\omega_1}^{\nu_1}\chi_{\omega_2}^{\nu_2}\dots\chi_{\omega_n}^{\nu_n}\,,
\quad\text{where}\quad
\nu_1m_1+\nu_2m_2+\cdots+\nu_nm_n\leq N\,,
\end{gather*}
where $N = M+1$ for $f$ and $N=M$ for $g$.

Since $\chi_{\omega_1}^{\nu_1}\chi_{\omega_2}^{\nu_2}\dots\chi_{\omega_n}^{\nu_n}$ decomposes into a linear combination of characters
$\sum a_\lambda
     \chi_{\lambda_1\omega_1+\cdots+\lambda_n\omega_n}$,
where $\lambda_j\leq\nu_j$ for all $j$, and $a_\nu=1$ (i.e. the exponent
$\nu= \nu_1\omega_1+\nu_2\omega_2+\cdots+\nu_n\omega_n$ appears in the sum with multiplicity $1$), we need only prove the Theorem for 
$f(X_1,\dots,X_n)= \chi_\lambda = \chi_{\lambda_1\omega_1+\cdots+\lambda_n\omega_n}$
and $g(X_1,\dots,X_n)= \chi_\mu = \chi_{\mu_1\omega_1+\cdots+\mu_n\omega_n}$.
Notice that it is quite possible for $\chi_0$, which is the constant function $1_\TT$,
to appear here.

\smallskip

We have from \eqref{wif} and \eqref{dfa}, 
\begin{eqnarray} \label{orthoRel}
\delta_{\lambda,\mu} &=& 
\frac{1}{|W|} \int_\TT \chi_\lambda \overline{\chi_\mu} \, |S_\rho|^2 \,d\theta_\TT 
= \frac{1}{|W|}\int_\TT S_{\lambda+\rho} \overline{S_{\mu+\rho}}\, d\theta_{\TT}
= \int_F S_{\lambda+\rho} \overline{S_{\mu+\rho}}\, d\theta_{\TT} \\
&=& \frac1{c_G(M+h)^n}\sum_{x\in {F^\circ_{M+h}}} S_{\lambda+\rho}(x) \overline{S_{\mu+\rho}(x)}\,,
\nonumber
\end{eqnarray}
(which is exactly what we have to prove) as long as when $\lambda \ne \mu$ there are no pairs $w,v\in W$ for which 
$w(\lambda +\rho)(x) -v(\mu+\rho)(x) \in \ZZ$ for all $x\in \frac{1}{M+h} \check P$, or, as
we pointed out above, 
$$
w(\lambda +\rho) -v(\mu+\rho) \in (M+h)\,Q\,.
$$
We will now show this cannot happen.

Consider the weights $w(\lambda +\rho)$, $w\in W$. These weights are all of the form
$\lambda + \rho - \beta$ where $\beta$ is a sum of positive roots (including the case 
when it is the empty sum, $0$). Now $-\check\a$ is the highest co-root and so is actually
a dominant co-weight. This gives us that $\l \beta,-\check\a\r \ge 0$ and so 
$$
m\mbox{-deg}\,(w(\lambda +\rho)) \le m\mbox{-deg}\,(\lambda +\rho) \le M+h\,.
$$
The lowest weight in the $W$-orbit of $\lambda + \rho$ is $w_{opp}(\lambda + \rho)$
and its $m$-degree is the negative of the $m$-degree of $\lambda + \rho$ as we
saw above. All the other weights in its orbit are of the form $w_{opp}(\lambda + \rho) +\beta$
for some sum $\beta$ of positive roots and this then gives us 
$$
-(M+h) \le -m\mbox{-deg}\,(\lambda +\rho) \le m\mbox{-deg}\,(w(\lambda +\rho)) .
$$
In short
$$
 -(M+h) \le -m\mbox{-deg}\,(\lambda +\rho)\le m\mbox{-deg}\,(w(\lambda +\rho))
\le  m\mbox{-deg}\,(\lambda +\rho)\le M+h \,.
$$
Exactly the same holds for $v(\mu+\rho)$ except that the inequalities are now
strict since the degree of $g$ is at most $M$. Combining, we obtain
$$
 -2 (M+h) < w(\lambda +\rho) -v(\mu+\rho) < 2(M+h)
 $$
for all $w,v \in W$. 

Now for any choice of $w,v$, each element in the $W$-orbit of
$w(\lambda +\rho) -v(\mu+\rho)$ is another element of the same form, and so its
$m$-degree is constrained in the same way. Thus if there is a pair $w,v$ for which
$w(\lambda +\rho) -v(\mu+\rho) \in (M+h)\,Q$ then, since $(M+h)Q$ is
$W$-invariant, we can assume that 
$\phi:=  w(\lambda +\rho) -v(\mu+\rho) \in P^+ \cap (M+h)\,Q$, i.e. it is dominant.
But then if $\phi \ne0$ it contradicts Lemma \ref{keyLemma}. It follows that $\phi=0$,
and then due to the fact that $\lambda+\rho$ and $\mu+\rho$ are strictly dominant,
we have $\lambda = \mu$. 

This proves the separation condition holds, and finishes
the proof of the Theorem.
 \end{proof}

\subsection{Duality}\ 

The essence of the cubature formula of Corollary~\ref{cubature} is the duality between
dominant weights of $m$-degree not exceeding $M$ and the elements of finite order
$M+h$ arising from the fundamental region $F$. If we use the fact that
the $m$-degree of $\rho$ is $h-1$ and the fact that any regular EFO in $F$
can be expressed as $\check \mu + \check\rho$ (or $x_{\check\mu+\check\rho}$
in $\TT$), where $\check \mu$
is co-dominant, then the cubature matrix is :
$$
\left(M_{\lambda\ \check\mu}\right)=\left( S_{\lambda+\rho}(x_{\check\mu + \check \rho}) \right) \, ,
$$
where $\lambda,\check \mu$ run over all solutions to the equations
\begin{eqnarray*}
\sum \check m_j (\lambda_j+1) &=& M+h\\
\sum  m_j (\check\mu_j+1) &=& M+h \, .
\end{eqnarray*}
Since the marks and co-marks are simply permutations of each other, the
symmetry of this pairing is completely manifest: if we order the co-ordinates
on each side to take account of this permutation, the solutions to the
two equations look identical and the matrix becomes symmetric. Moreover, formally
the situation is the same for $G$ and the corresponding group $\check G$, with the roles
of characters and elements of finite order interchanged.


\section{Approximating functions on $\Omega$}\label{Approximations}

To simplify notation and the following discussion we introduce an inner product on the space
$L^2_K(\Omega)$ of all complex-valued functions $f$ on $\Omega$ for which 
$\int_\Omega |f|^2 \,K^{1/2} < \infty$. It is defined by 
\begin{equation}\label{innerProduct}
\langle f, g \rangle_K: = \int_\Omega f \overline g \, K^{1/2} =  
\int_\Omega f(X_1,\dots,X_n) \overline{g(X_1,\dots,X_n)} \, K^{1/2}(X_1,\dots,X_n) \, dX_1\dots dX_n \,.
\end{equation}

{\em Note here that the definition of conjugation is made in terms of \,{\rm Remark~\ref{conjugation!}}.}

\smallskip
Since $K^{1/2}$ is continuous and strictly positive on $\Omega$, $\langle f, f\rangle_K \ge 0$
with equality if and only if $f$ is zero almost everywhere in  $\Omega$ (relative to 
Lebesgue measure). Relative to this inner product $L^2_K(\Omega)$ is a Hilbert
space.  We call $\langle f,f\rangle_K^{1/2}$ the $L^2_K$-norm of $f$. 

The results of \eqref{orthoRel} 
show that the polynomials defined by 
$$
X_\lambda := \chi_\lambda(x), \quad x\in F^\circ \, ,
$$
where $\lambda = \sum \lambda_j \omega_j $ runs through the set of dominant weights $P^+$, form an orthonormal set:
$$
 \langle X_\lambda, X_\mu \rangle_K = \delta_{\lambda, \mu} \, .
 $$
These polynomials are just those that result from expanding $\chi_\lambda$
as a polynomial (of degree equal to the $m$-degree $|\lambda|_m$ of $\lambda$) in the fundamental characters $X_j = \chi_{\omega_j}$.

The set $\{X_\lambda\}$ of actually forms a Hilbert basis for $L^2_K(\Omega)$ (the main
point being that they actually span the entire space). This can be
seen by relating functions $f$ on $\Omega$ back to functions $\tilde f$ on $F^\circ$
through $\tilde f(x) :=   f(X_1, \dots, X_n) = f(\chi_{\omega_1}(x), \dots, \chi_{\omega_n}(x))$.
Then 
\begin{equation}\label{int}
\infty > \int_\Omega |f|^2 K^{1/2} = \int_{F^\circ} |\tilde f|^2(x) \, S_\rho(x) \overline{S_\rho(x)} \, dx\\
 \, .
 \end{equation}
Formally $\tilde f$ exists only on $F^\circ$, but
we can extend it to a function on all of $\TT$ by $W$-symmetry if necessary. This might
make it a bit easier to relate to  \S\ref{Winvariance}, which we are now going to use.

We have already seen the integral on the right hand side of \eqref{int} in \eqref{newInnerProduct},
and according to \eqref{FEoff} and \eqref{rewriteInnerProduct} we have 
$$ \tilde f = \sum_\lambda b_\lambda\, \chi_\lambda = \sum_\lambda b_\lambda\, \widetilde{X_\lambda} \,,
$$
with 
$$
b_\lambda = \frac{1}{|W|}\int_\TT \tilde f \,\overline{\chi_\lambda}\, |S_\rho|^2 \, d\theta_\TT 
= \int_{F^\circ} \tilde f\, \overline{\chi_\lambda} \,|S_\rho|^2 \, d\theta_\TT
= \int_\Omega f \overline{X_\lambda} \,K^{1/2} =\langle f, X_\lambda \rangle_K\, ,
$$ 
(see \eqref{integralsOverT}).
 \smallskip

Thus for each $f \in L^2_K(\Omega)$ we have its expansion
$$
 f \bumpeq \sum_\lambda \langle f, X_\lambda\rangle_K  \, X_\lambda \, ,
 $$
where the $\bumpeq$ means equality in the sense of equality Lebesgue-almost-everywhere. 

The truncated sums
$$
\sum_{|\lambda|_m \le M} \langle f, X_\lambda\rangle_K  \, X_\lambda \, ,
$$
where $|\lambda|_m$ stands for the $m$-degree of $\lambda$,
are polynomials of $m$-degree at most $M$ in the variables $X_1, \dots, X_n$.

\subsection{Optimality of the approximation}\ 

We now show that these polynomials are the best possible approximations for $f$
in terms of the $L^2_K$-norm. This result
is in fact a natural consequence that one would expect from the situation that 
we have created here, since it essentially relates back to the Fourier analysis
of $\TT$. See also \cite{Xu}. 

\begin{prop}\label{approxThm}\ 

 Let $f \in L_K^2(\Omega)$. Amongst all polynomials $p(X_1,\dots,X_n)$ of $m$-degree less
than or equal to $M$, the polynomial $\sum_{|\lambda |_m \le M} \langle f, X_\lambda\rangle_K  \, X_\lambda $ is the best approximation to $f$ relative to the $L^2_K$-norm.
\end{prop}
 \begin{proof}
 Let $p = \sum_{|\lambda|_m \le M}  b_\lambda \, X_\lambda $ be any polynomial.
Let $a_\lambda:= \langle f, X_\lambda\rangle_K$ for all $\lambda$ and set
$q:= \sum_{|\lambda|_m \le M}  a_\lambda \, X_\lambda $. 
 \begin{eqnarray*}
 \langle f -p, f-p\rangle_K &=& \langle f , f\rangle_K - \langle f , p\rangle_K - \langle p, f\rangle_K
 +\langle p, p\rangle_K\\
 &=& \langle f , f\rangle_K - \sum_\lambda a_\lambda \overline{b_\lambda} -
  \sum_\lambda b_\lambda \overline{a_\lambda} + \sum_\lambda |b_\lambda|^2\\
  &=&  \langle f , f\rangle_K - \sum_\lambda |a_\lambda|^2 + 
  \sum_\lambda |a_\lambda|^2- \sum_\lambda a_\lambda \overline{b_\lambda} -
  \sum_\lambda b_\lambda \overline{a_\lambda} + \sum_\lambda |b_\lambda|^2\\
  &=& \langle f -q , f -q\rangle_K   +\sum_\lambda |b_\lambda - a_\lambda|^2 
  \ge  \langle f -q , f -q\rangle_K\,
 \end{eqnarray*}
 with equality if and only if $b_\lambda = a_\lambda$ for all $\lambda$. 
 \end{proof}

\section{Example: A cubature formula for $G_2$}\label{G2}
We consider an example where the Lie group is the exceptional simple Lie group of type $G_2$ and the value of $M$ is chosen to be $M=8$. 
Before we really start, let us recall pertinent information about $G_2$ and related objects. 
Then we can start to solve the problem of our example. There are three steps to it. First the $S$-functions of $m$-degree $M+1=9$ are to be found (or more specifically the highest
weights that label them). Then their common zeros of Ad-order $M+h=14$ in the interior of $F$ are to be determined, and then only the cubature formula (matrix) can be written down.

\subsection{Pertinent data about $G_2$}\ 

From the $G_2$ diagram in Figure~\ref{diagrams} we find the following:
\begin{alignat*}{2}
&\l\alpha_1,\alpha_1\r\ :\ \l\alpha_2,\alpha_2\r=3 : 1\,,&\quad
   &A=\left(\begin{smallmatrix}2&-3\\-1&2\end{smallmatrix}\right)\,;\\
&\l\check\alpha_1,\check\alpha_1\r\ :\ \l\check\alpha_2,\check\alpha_2\r=1:3\,,&\quad
   &\check A=\left(\begin{smallmatrix}2&-1\\-3&2\end{smallmatrix}\right)\,;\\
  -&\alpha_0 = 2\alpha_1+3\alpha_2\,,&\qquad 
    -& \check\alpha_0= 3\check\alpha_1+2\check\alpha_2\,;\\
&\mbox{marks}\quad {2,\ 3}\,,&\qquad &\mbox{co-marks}\quad {3,\ 2}\,;\\
&\mbox{Coxeter number}\; h=6\,.&&
\end{alignat*}

Therefore we have
\begin{alignat*}{4}
\alpha_1&=2\omega_1-3\omega_2\,,&\qquad 
\alpha_2&=-\omega_1+2\omega_2\,,&\qquad
\omega_1&=2\alpha_1+3\alpha_2\,,&\qquad
\omega_2&= \alpha_1+2\alpha_2\,,\\
\check\alpha_1&=2\check\omega_1-\check\omega_2\,,&\qquad 
\check\alpha_2&=-3\check\omega_1+2\check\omega_2\,,&\qquad
\check\omega_1&=2\check\alpha_1+\check\alpha_2\,,&\qquad
\check\omega_2&=3\check\alpha_1+2\check\alpha_2\,.
\end{alignat*}
The link between the bases is the $\ZZ$-duality requirement 
$$
\l\alpha_j,\check\omega_k\r=\l\check\alpha_j,\omega_k\r=\delta_{jk}\,,\qquad j,k=1,2\,.
$$
The set of positive roots $\Delta_+$ and their half sum $\rho$ are:
$$
\Delta_+=\{\alpha_1,\ \alpha_2,\ \alpha_1+\alpha_2,\ \alpha_1+2\alpha_2,\ 
        \alpha_1+3\alpha_2,\ 2\alpha_1+3\alpha_2\},\quad
\rho=3\alpha_1+5\alpha_2=\omega_1+\omega_2\,.
$$
The fundamental domain of $G_2$ is the convex hull of its three vertices 
$\{0,\tfrac12\check\omega_1,\tfrac13\check\omega_2\}$.

\subsection{Finding $S$-functions of $m$-degree 9}\ 

The $m$-degree of a $G_2$ weight $\lambda= \{\lambda_1,\lambda_2\} =\lambda_1\omega_1+\lambda_2\omega_2$ is calculated as $\check m_1\lambda_1+\check m_2\lambda_2=3\lambda_1+2\lambda_2$.
We find first all the weights with $m$-degree $\leq9$. They are the following,
\begin{align*}
\{0,0\},\,\{0,1\},\,\{0,2\},\,\{0,3\},\,\{0,4\},\,\{1,0\},\,
\{1,1\},\,\{1,2\},\,\{1,3\},\,\{2,0\},\,\{2,1\},\,\{3,0\} \,.
\end{align*}
Among these there are just two weights with $m$-degree equal to 9, namely $\{1,3\}$ and $\{3,0\}$. Thus there are only two $S$-functions we need to consider, $S_{1,3}(x)$ and $S_{3,0}(x)$. 

\subsection{Finding common zeros of the $S$-functions of $m$-degree 9}\ 
 
The interpolation points are the points of the set $F^\circ_{14}$ of EFO's of $Ad$-order $M+h=14$ which are in the interior of $F$. They must satisfy the sum rule $14=s_0+2s_1 +3s_2$. We write them in the in Kac coordinates $[s_0,s_1,s_2]$  as well as in the $\check\omega$-basis: 
\begin{equation}\label{zeros}
\begin{alignedat}{4}
[9,1,1]&=(\tfrac1{14},\tfrac1{14}),\quad &
[7,2,1]&=(\tfrac1{7},\tfrac1{14}),\quad &
[5,3,1]&=(\tfrac3{14},\tfrac1{14}),\quad &
[3,4,1]&=(\tfrac27,\tfrac1{14}),\\
[1,5,1]&=(\tfrac5{14},\tfrac1{14}),\quad &
[6,1,2]&=(\tfrac1{14},\tfrac17),\quad &
[2,1,1]&=(\tfrac17,\tfrac17),\quad &
[2,3,2]&=(\tfrac3{14},\tfrac17),\\
[3,1,3]&=(\tfrac1{14},\tfrac3{14}),\quad &
[1,2,3]&=(\tfrac1{14},\tfrac3{14})
\end{alignedat}
\end{equation}
The strict adjoint order of all the EFO's but one is $14$. It is $7$ for 
$\tfrac17\check\omega_1+\tfrac17\check\omega_2$. Other EFO's of adjoint order $14$ are not shown in \eqref{zeros} because they are on the boundary of $F$. The boundary points are easily discarded because at least one of their Kac coordinates has to be zero. At these points the $S$-functions 
$S_{1,3}$ and $S_{3,0}$ simultaneously vanish.

\subsection{Cubature formula}\ 

Underlying the cubature formula (for $M=8$ in $G_2$) is the fact that the matrix 
$$
X =X^{(8)}=\left( S_{\lambda + \rho}( x_{\check \mu}) \right)
$$
where the $\lambda$ run over the $10$ weights of $m$-degree at most $8$
and the $\check\mu$ run over the $10$ regular EFOs of $Ad$-order $14$, satisfies
$\frac{1}{(8+6)^2} X\,X^T = I_{10\times 10}$. Here is the matrix $X^{(8)}$
with entires rounded to $4$ figures. 

$$
\left(
\begin{array}{cccccccccc}
 -0.604 & -2.714 & -5.494 & -6.098 &
   -2.714 & -3.384 & -7. & -6.098 &
   -3.384 & -2.11 \\
 -2.714 & -7.604 & -6.098 & 1.506 &
   2.714 & -4.89 & 0. & 6.098 & 3.384
   & 3.384 \\
 -5.494 & -6.098 & 2.11 & -3.384 &
   -6.098 & 2.714 & 7. & -3.384 &
   2.714 & -0.604 \\
 -6.098 & 1.506 & -3.384 & -4.89 &
   6.098 & 7.604 & 0. & 3.384 &
   -2.714 & -2.714 \\
 -2.714 & 2.714 & -6.098 & 6.098 &
   -7.604 & 3.384 & 0. & 1.506 &
   -4.89 & 3.384 \\
 -3.384 & -4.89 & 2.714 & 7.604 &
   3.384 & -1.506 & 0. & -2.714 &
   -6.098 & -6.098 \\
 -7. & 0. & 7. & 0. & 0. & 0. & -7. &
   0. & 0. & 7. \\
 -6.098 & 6.098 & -3.384 & 3.384 &
   1.506 & -2.714 & 0. & -4.89 &
   7.604 & -2.714 \\
 -3.384 & 3.384 & 2.714 & -2.714 &
   -4.89 & -6.098 & 0. & 7.604 &
   -1.506 & -6.098 \\
 -2.11 & 3.384 & -0.604 & -2.714 &
   3.384 & -6.098 & 7. & -2.714 &
   -6.098 & 5.494
\end{array}
\right)
$$
Direct computation shows that $X^{(8)}(X^{(8)})^T=14^2I_{10\times 10}$.

\smallskip
The values of the weighting $K$ function at the interpolation points are given by

\medskip
\begin{center}
\small{
$0.364666,\ 7.36467,\ 30.1836,\ 37.1836,\ 7.36467,\ 
 11.4517,\ 49.,\ 37.1836,\ 11.4517,\ 4.45175$
}
\end{center}
\medskip

The image of the fundamental domain $F$ under the  mapping 
$x\mapsto (X_1(x), X_2(x))$ in $\RR^2$ along with the images of the 
$8$ EFOs that we have used in this example is shown in Fig.~\ref{omega}:
\begin{figure}
\centering
\includegraphics[width=12cm]{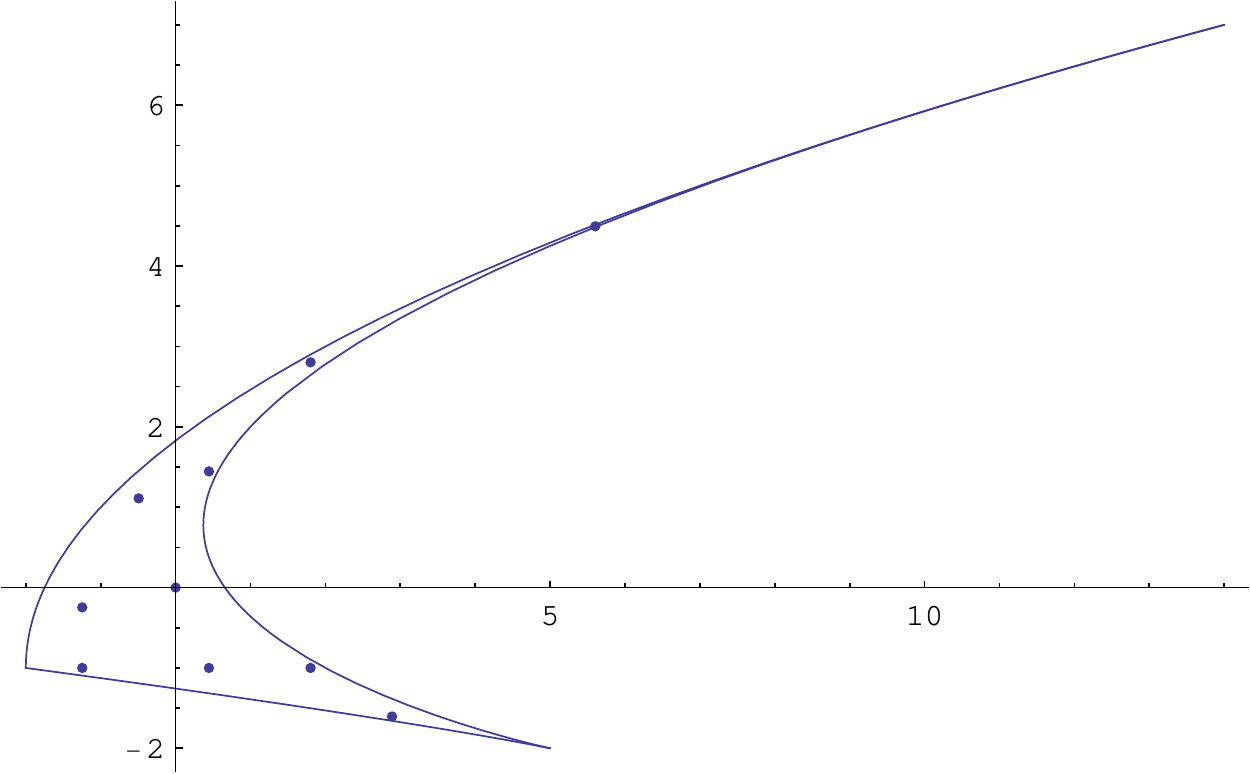}
\caption{The region $\Omega$ along with the $10$ regular EFOs of the example.}
\label{omega}
\end{figure}
The way in which EFOs fill out the region $\Omega$ is made clearer by the 
Fig.~\ref{filledOmega} which shows the distribution of the EFOs of Ad-order 106:
\begin{figure}
\centering
\includegraphics[width=12cm]{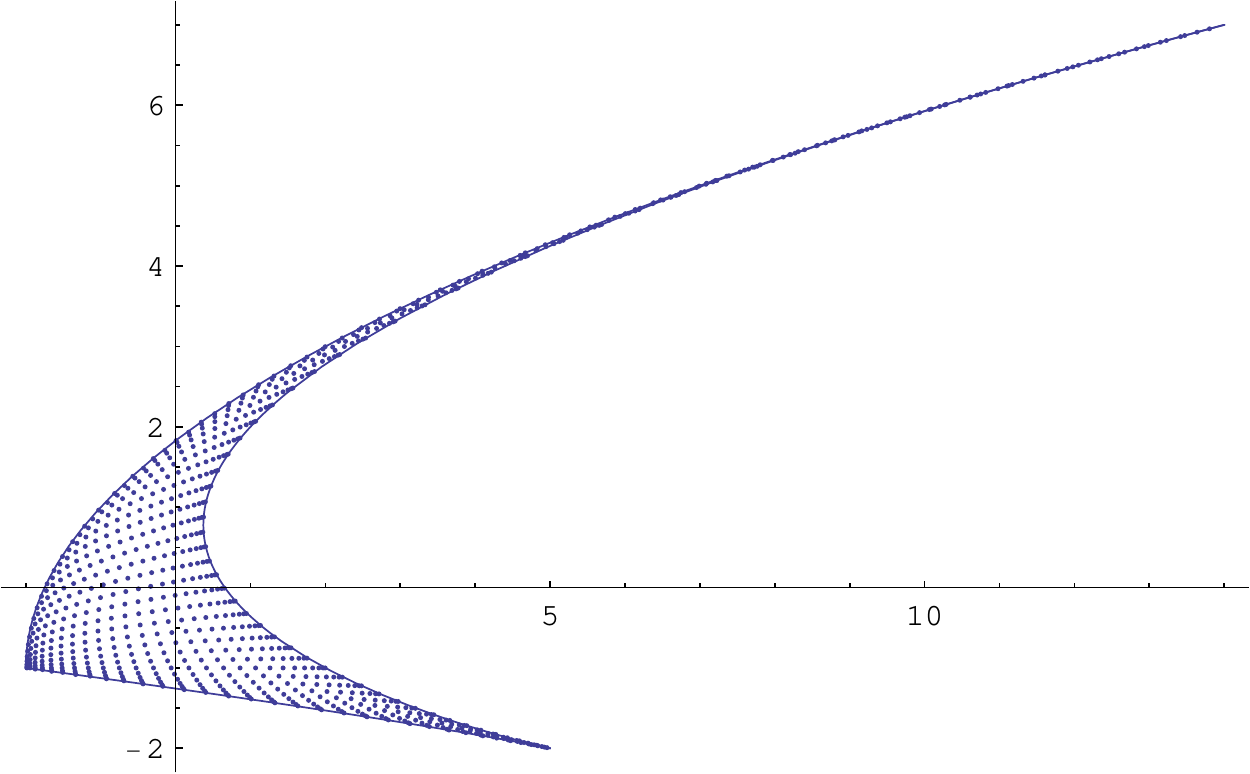}
\caption{The region $\Omega$ along with the $884$ regular EFOs of Ad-order $106$.}
\label{filledOmega}
\end{figure}

Examination of the graphs of the two functions $S_{1,3}$ and $S_{3,0}$ reveals
that they have (at least) $2$ other common zeros in $F$. These do not seem to be
related to EFOs. 

\subsection*{Acknowledgments}\ 

 RVM and JP gratefully acknowledge
the support of this research by the Natural Sciences and Engineering Research
Council of Canada and by the MIND Research Institute of Santa Ana, California.
We also wish to thank Prof. Yuan Xu who sent us the preprint \cite{LX} that 
made us aware of the potential connection between cubature formulae and 
root systems. 


\end{document}